\documentclass[a4paper]{article}
\usepackage{a4wide}

\usepackage{endnotes}
\usepackage{tikz}
\let\footnote=\endnote

\usepackage{bbm,amsmath,amsfonts,amssymb,amsthm}
\usepackage{latexsym}
\usepackage{epsfig, psfrag}
\usepackage{dsfont}
\usepackage[ruled,vlined]{algorithm2e}
\usepackage[font=small,labelfont=bf,labelsep=none]{caption}
\usepackage{stmaryrd}
\usepackage{figlatex}
\usepackage{graphics}
\usepackage{multicol}

\graphicspath{{Figures/}}

\newcommand{\Pb}{{\mathds{P}}}
\newcommand{\Eb}{{\mathds{E}}}

\newcommand{\cE}{\mathcal{E}}
\newcommand{\cS}{\mathcal{N}}
\newcommand{\cX}{\mathcal{X}}

\newcommand{\ft}{f_T} 
\newcommand{\fp}{f_P} 
\newcommand{\Ft}{F_T} 

\newcommand{\ffp}{g} 
\newcommand{\overffp}{\ffp^{\sup}} 
\newcommand{\underffp}{\ffp^{\inf}} 
\newcommand{\FFpp}{H} 
\newcommand{\FFppinf}{\underline{\FFpp}} 
\newcommand{\FFppsup}{\overline{\FFpp}} 
\newcommand{\Inf}{h^{\inf}} 

\newcommand{\proj}{\psi}

\newcommand{\Ysup}{Y^{\sup}} 
\newcommand{\Yinf}{Y^{\inf}} 

\newcommand{\N}{\mathbb{N}}
\newcommand{\integer}{\mathbb{N}}
\newcommand{\R}{\mathbb{R}}

\def\ph{f}
\def\Ph{F}
\def\SS{\mathcal{X}}
\def\Space{\mathcal{X}}
\def\evtSet{\mathcal{E}}
\def\Evt{\mathcal{E}}
\def\A{\mathcal{E}}
\def\M{\mathcal{M}}
\def\init{\nu}
\newcommand{\bydef}{\stackrel{\rm{def}}{=}}
\def\peq{\preceq}
\def\I{\mathcal I}
\def\J{\mathcal J}

\newcommand{\ie}{{\it i.e.~}}
\newcommand{\iid}{{\em i.i.d.~}}
\newcommand{\Iid}{{\em I.i.d.~}}

\newcommand{\PcX}{\mathcal{P(\cX)}}
\newcommand{\Ux}{S_{x}}

\newcommand{\st}{\preceq_{st}}

\usepackage{url}
\graphicspath{{Figures/}}

\newtheorem{theorem}{Theorem}
\newtheorem{lemma}[theorem]{Lemma}
\newtheorem{proposition}[theorem]{Proposition}

\newtheorem{definition}[theorem]{Definition}
\newtheorem{example}[theorem]{Example}

\newtheorem{remark}[theorem]{Remark}

\begin{document}




\title{Exact Simulation for Assemble-To-Order Systems}

\author{Ana {\sc Bu\v si\'c}\thanks{INRIA and Computer Science Department of \'Ecole Normale Sup\'erieure (DI ENS), Paris, France. 
E-mail: {\tt ana.busic@inria.fr}.}  \and
Emilie {\sc Coupechoux} \thanks{Laboratoire I3S, Universit\'e Nice Sophia Antipolis, Paris, France. 
E-mail: {\tt coupecho@i3s.unice.fr}. This paper is part of the author's PhD thesis done at INRIA-ENS.}
}

\maketitle

\begin{abstract} 
We develop exact simulation (also known as perfect sampling) algorithms  for a family of assemble-to-order systems. Due to the finite capacity, and coupling in demands and replenishments, known solving techniques are inefficient for larger problem instances. 
We first consider the case with individual replenishments of items, and derive an event based representation of the Markov chain that allows applying existing exact simulation techniques, 
using the monotonicity properties or bounding chains. 
In the case of joint replenishments, the state space becomes intractable for the existing methods. 
We propose new exact simulation algorithms, based on aggregation and bounding chains, that allow a significant reduction of the state space of the Markov chain. 
We also discuss the coupling times of considered models and provide sufficient conditions for linear (in the single server replenishment case) or quadratic (many server case) complexity 
of our algorithms in terms of the total capacity in the system.
\end{abstract} 




%


\section{Introduction}

Exact simulation (also called perfect sampling) algorithms
draw unbiased samples from a target distribution. This distribution is in general not known, or cannot be efficiently computed due to:
\begin{itemize}
 \item the computational complexity of the normalizing constant
of an otherwise known distribution (e.g. in statistical physics applications, 
approximation algorithms of $\sharp P$-complete problems, or analysis of queueing systems with product form solution), or 
\item the target distribution that is a limiting distribution of a Markov chain that cannot be efficiently solved by analytic methods 
(e.g. queueing systems not having a product form solution). 
\end{itemize}
We will focus mainly on the latter case and develop exact simulation algorithms for Assemble-To-Order (ATO) systems.  

\cite{propp96exact} used a \emph{coupling from the 
past} scheme to derive an exact simulation algorithm - called PSA (Perfect Sampling
Algorithm) in the following - providing unbiased samples from the stationary distribution of an ergodic 
Markov chain with a finite state space. 
Many variants of their algorithm have been developed since in
various contexts. We give in Section \ref{sec:PSA} a brief overview of PSA and mention some works directly linked to the present article (for more information, see the annotated bibliography by \cite{biblio}). 

PSA provides unbiased samples from the stationary distribution of an ergodic Markov chain 
in finite expected time. This is its main advantage over Markov Chain Monte Carlo (MCMC) simulation methods that construct one trajectory of the chain and stop after some long enough burn in period, when the distribution 
of the current 
state 
is estimated to be close enough to the stationary distribution 
(see \cite{AsmGly07} for an overview of MCMC methods). 
This burn in period depends on the mixing time of the chain, that is in general difficult to compute.  
The stopping criterion for MCMC is thus often based on pessimistic bounds, that can be much larger than the coalescence time of the coupling from the past scheme used by PSA 
(note that the coalescence time provides an upper bound for the mixing time, see for instance \cite{Peres}). 
Moreover, PSA detects the exact coalescence time; there is no explicit
need for its estimation for the exactness of the algorithm. This represents a clear advantage over MCMC methods using
stopping criteria based on the mixing time, when the latter is not known.


The efficiency of PSA \textit{a priori} depends on the size of the state space. When the system is monotone, one can easily overcome this issue (as mentioned in \cite{propp96exact} 
and explained in Section \ref{sec:PSA}). 
Developing efficient exact simulation algorithms for non-monotone Markov chains with a very large state space remains challenging: 
Bounding chains can be used to detect coalescence, as in \cite{KendallMoeller-2000}, \cite{Hu04} or \cite{ashes}.

\cite{KendallMoeller-2000} describe the general idea and provide the construction of bounding chains for spatial birth-and-death processes. \cite{Hu04} gives a general approach for Markov chains with local interactions, including 
Gibbs or Metropolis-Hastings samplers.
Envelope Perfect Sampling Algorithm (EPSA) developed by \cite{ashes} 
gives a rather general construction of bounding chains, under the assumption that the state space is a lattice. 

ATO systems with joint returns, that motivated our work, have a very large state 
space (that grows exponentially with the number of item types). Furthermore the state space is \textit{not} naturally equipped with a lattice order relation, so  EPSA cannot be used. 
Our goal is two-fold: 
\begin{itemize}
 \item We propose a new exact simulation method, generalizing EPSA and based on aggregation. This Aggregated Envelope Perfect Sampling Algorithm (AEPSA) can be applied to reduce 
the state space and overcome the lack of a lattice structure.
 \item We provide a detailed treatment of exact simulation algorithms for ATO systems: 
More precisely, we use known algorithms (PSA or EPSA) for ATO systems with individual returns, and our new algorithm (AEPSA) for ATO systems with joint returns. Each type of ATO system considered is described below, 
and the choice of the algorithm (PSA, EPSA or AEPSA) depending on the system is summarized in Table \ref{table:A-E-PSA}.
\end{itemize}

Up to our knowledge, this is the first time perfect sampling techniques are applied to ATO systems 
(for an overview on ATO systems, see \cite{SZ03}).
We focus here on {\it continuous-review models}, with 
{\it exponential replenishment times} and {\it finite stock capacities}. 
As is common in the ATO literature, we assume at most one component of each type will be demanded for any item (the unit demand case). 
We consider two different options for the out-of-stock situation: A demand can be fulfilled partly (just the components that are available), referred to as a partial order service (POS); 
or lost fully, referred to as a total order service (TOS).  
Also, 
we distinguish between two different situations for the replenishment/return of components: components are either returned individually 
or jointly.
We briefly mention only the possible 
solution techniques that are directly related to the models we consider in this paper (an overview is given in Table \ref{table:framework}).

\begin{table}[htb]
\begin{center}
\begin{tabular}{|l|l|l|}
\hline
Returns & 1. TOS & 2. POS  \\
\hline
A. Individual & {\bf Exact} (matrix geometric):  & {\bf Exact} (matrix geometric):\\
 & \cite{SXL99} & \cite{SXL99}  \\
 &  & {\bf Bounds}: \\
 & & \cite{Xu99}, \cite{LX00}, \cite{Xu02}, \\
 &  & \cite{DSX03} \\
\hline
B. Joint & {\bf Exact} (product form): & {\bf Approximations}:  \\ 
& \cite{Kel91} & \cite{VH09}  \\
 &  & {\bf Bounds}:  \\
 &  & \cite{BVSW12} \\
\hline\end{tabular}
\caption{\label{table:framework} Solution methods for ATO systems.}
\end{center}
\end{table}

\cite{SXL99} proposed an exact evaluation of ATO systems with individual returns, both for TOS and POS, by using a matrix geometric approach. This exact method, however, is computationally inefficient for larger problem instances. 

\cite{Xu99}, \cite{LX00}, \cite{Xu02} studied the effect of correlation (for the arrival process) on a variety of system performance measures for correlated queueing systems, including ATO-POS systems. 
\cite{DSX03} presented several approximations and bounds on the performance of ATO-POS systems with individual returns. 

The best known example of models that have joint returns of resources and TOS are telecommunication systems, or specifically loss networks (\cite{Kel91}). 
In these networks, demands arrive, for example a phone call, that need several links to be simultaneously available. If all links are available, the call is connected. After the call is finished, all links are simultaneously released. When one or more of the links is not available, the call does not connect, and the demand for all links is lost. 
Although loss networks have a product-form solution, exactly computing the blocking probabilities for this system
is known to be a difficult problem (\cite{LMK93}), due to the normalizing constant.

An example of ATO systems with joint returns and a partial order service is the service tool problem, considered in \cite{VH09}. 
In this problem, to perform a maintenance action, several service tools are needed at the same time. After usage, all tools return to the location they were sent from together. 
Whenever one or more tools are not present, they are sent by an emergency shipment to enable the initiation of the maintenance action as soon as possible. 
For the supply location under consideration the demand for these emergency shipped tools is lost. 
\cite{VH09} developed different approximations. Some of these approximations provide provable bounds (\cite{BVSW12}). 
For larger instances, however, these bounds are still time consuming.

\begin{table}[htb]
\begin{center}
\begin{tabular}{|l|c|c|}
\hline
Returns & 1. TOS & 2. POS  \\
\hline
A. Individual & Envelope Perfect Sampling Algorithm: & Perfect Sampling Algorithm: \\
 &  \cite{BGV08,ashes} & \cite{propp96exact} \\
\hline
B. Joint & \multicolumn{2}{c|}{Aggregated Envelope Perfect Sampling Algorithm:} \\
& \multicolumn{2}{c|}{Section \ref{sec:method}} \\
\hline\end{tabular}
\caption{\label{table:A-E-PSA} Perfect sampling algorithms for ATO systems.}
\end{center}
\end{table}
 
For larger instances, all four cases are difficult to analyze directly and we will discuss in this paper how to develop exact simulation algorithms for ATO systems, as an 
alternative for bounding techniques developed in the literature. In Table~\ref{table:A-E-PSA}, we mention which algorithm (PSA, EPSA or AEPSA) is used for each of the four models. The monotonicity of the POS system with individual returns (Proposition~\ref{prop:ATOpos_monotone}) allows the use of PSA, while we apply EPSA to the TOS system with individual returns, which is non-monotone (Proposition~\ref{prop:I3}). To handle the case of joint returns, we use the method developed in Section~\ref{sec:method} (AEPSA). 
Loss networks (case B1) have a product-form solution, which makes the exact calculation easier to some extent. We therefore focus on the other cases and we will only briefly mention 
in Section~\ref{sec:frc} how we can adapt the approach developed for the POS case with joint returns (Section~\ref{sec:Aggreg_POS}) to the TOS case. 

In addition, we give bounds on the complexity of our algorithms.
The complexity of the exact simulation algorithms we develop depends on the coupling time 
that is usually difficult to estimate and even to bound, 
except for some specific Markov chains 
(\cite{Peres}).
In the context of queueing networks, \cite{Dopper}
have shown that the coupling time in an M/M/1/C queue is linear in capacity $C$, when $\lambda \not=\mu$
(with $\lambda$ being the arrival and $\mu$ the service rate), and quadratic when $\lambda =\mu$.
They used this fact to derive an upper bound for the coupling time of an acyclic network of $K$ M/M/1/C queues that is $O(KC^2)$. These results have been extended to cyclic networks under some additional 
hyperstability conditions (\cite{AnselmiPS11}). 
However, these results use both the fact that the system is monotone under the usual product partial order and that each event in the system can only influence up to two different components. 
In our case, joint arrivals and services can modify many components at the same time. Also,
we do not always have the monotonicity property (see TOS case).

\begin{figure}

\begin{small}
\begin{center}
\begin{tabular}{cc}
 
\begin{tikzpicture}
  [scale=.8,auto=left]
  \node (n1) at (3,7)  {1};
  \node[draw,shape=circle] (n2) at (3,5)  {2};
  \node[draw,shape=rectangle] (n3) at (1,4)  {3};
  \node[draw,shape=rectangle] (n4) at (1,2)  {4};
  \node[draw,shape=circle] (n5) at (5,4)  {5};
  \node[draw,shape=rectangle] (n6) at (5,2)  {6};
  \node (n7) at (3,1)  {7};

  \foreach \from/\to in {n1/n2,n2/n3,n3/n4,n4/n7,n2/n5,n5/n6,n6/n7}
    \draw (\from) -- (\to);

\end{tikzpicture} &

\begin{tikzpicture}
  [scale=.8]
  \node (n1) at (3,7)  {\textcolor{white}{1}};
  \node (n2) at (3,5)  {\textcolor{white}{2}};
  \node (n3) at (1,4)  {\textcolor{white}{3}};
  \node (n4) at (1,2)  {\textcolor{white}{4}};
  \node[draw,shape=circle] (n5) at (3,4)  {\textcolor{white}{5}} ;
  \node[draw,shape=rectangle] (n6) at (3,2)  {\textcolor{white}{6}};
  \node (n7) at (3,1)  {\textcolor{white}{7}};
  \node[right]  at (4,4)%
     {Perfect sampling methods};
  \node[right]  at (4,2)%
     {Application to ATO systems};

\end{tikzpicture} \\

\end{tabular}
\caption{\label{fig:dependencies} Dependencies between sections.}
\end{center}\end{small}
\end{figure}
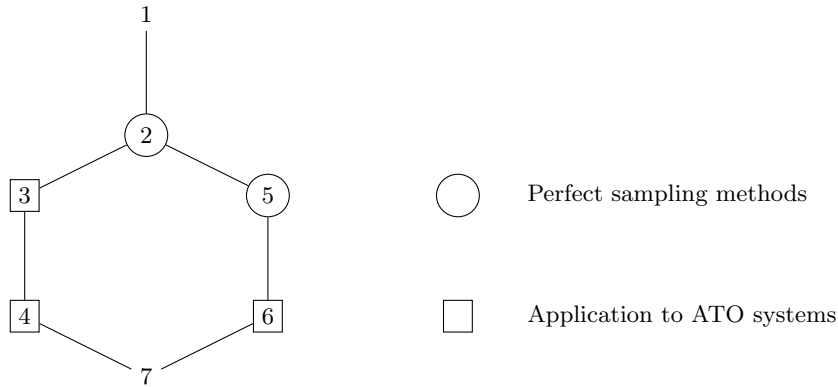

The paper is organized as follows. In Section \ref{sec:PSA} we give an overview of perfect sampling and the related literature. Sections \ref{sec:POS} and \ref{sec:TOS} 
are devoted respectively to ATO-POS and ATO-TOS systems with individual returns. 
In Section \ref{sec:method} we present AEPSA, that
we apply to the ATO-POS model with joint returns in Section \ref{sec:Aggreg_POS}. Finally, in Section \ref{sec:frc} we discuss some possible extensions of our work and provide conclusions. 
In Figure \ref{fig:dependencies}, we give the dependencies between sections.

\section{Perfect sampling and the method of envelopes}
\label{sec:PSA}

The evolution of a finite Discrete Time Markov Chain (DTMC) can always be obtained using a finite number of discrete 
events (or actions). We consider a system description similar to Generalized Semi Markov Processes (see \cite{cassandras}), 
with a focus on state changes rather than on time: we consider a tuple $\M = (\SS,\evtSet,\init,\ph)$ where $\SS$ is
a finite state space, $\A$ is the {\em set of events}, 
$\init$ is a probability distribution on $\A$, and $\ph$ is a {\em transition function}, 
$\ph : \SS \times \evtSet \to \SS$.

This transition function $\ph$ can be naturally extended to 
   words $a_{1\to t} \bydef a_1a_2\dots a_t \in \evtSet^t$, 
  $t \in \N$ (where $a_{1 \rightarrow 0} := \epsilon$ is the empty word). 
For any $t \in \N$, $\ph :{\SS \times  \evtSet^t} \to \SS$ is defined  by:
$\ph \left(x, \epsilon \right) \bydef x$ and $\ph \left(x, a_{1\to t} \right) \bydef\ph \left(\ph \left(x, a_{1\to t-1} \right), a_{t} \right)$ for $t\geq 1$.

Let $(a_t)_{t\geq 1}$ be an infinite \iid sequence of random events in $\evtSet$, distributed according to $\init$.
Then for any $x_0 \in \SS$, the random process 
$(X_t \bydef \ph(x_0,  a_{1 \rightarrow t}))_{t \geq 0}$
is a Markov chain started in $x_0$ with probability transition matrix $P$ given by:
\begin{equation}\label{eq:markov_automaton}
  \text{for all } x,y \text{ in } \SS, ~~ P(x,y) = 
{\sum_{a \in \evtSet, \; \ph(x,a) = y}} \init(a) .
\end{equation}
We say that the Markov chain $(X_t)$ is \textit{generated} by $\M$ and $(a_t)_{t\geq 1}$.

Conversely, for any probability transition matrix $P$ on a finite state space $\SS$,
it is easy to see that there exists a tuple $\M=(\SS,\evtSet,\init,\ph)$ such that 
\eqref{eq:markov_automaton} holds, \ie such that $\M$ generates a Markov chain
on $\SS$ with 
transition matrix $P$, but that representation is in general not unique.
However, such a representation
naturally arises for many systems, including Markovian queueing networks.

We can build a family of Markov chains 
$\{(X_t(x) = \ph(x, a_{1 \rightarrow t}))_{t \ge 0} \mid x \in \SS\}$ 
starting from each state $x \in \SS$, referred to as the {\em grand coupling} 
generated by $\M$ and $(a_t)_{t\geq 1}$ (\cite{Peres}).
We will say that the grand coupling has {\em coupled} (or more precisely {\em coalesced}) at time $t$
if all the Markov chains of the family have reached the same state. Using the notation $\ph \left(U, a_{1\to t} \right) \bydef \left\{ \ph \left(x, a_{1\to t} \right), x\in U \right\}$ for any subset $U\subset \SS$, this is equivalent to the fact that $\ph(\SS, a_{1 \rightarrow t})$ is reduced to a singleton. 
In the following, $|V|$ denotes the cardinality of set~$V$.

\subsection{Perfect Sampling} 

Let $\left(X_t\right)_{t \in \N}$ be an irreducible
and aperiodic DTMC with finite state space
$\SS$ and transition matrix $P$. 
Consider a discrete event system representation $\M=(\SS,\evtSet,\init,\ph)$ that satisfies 
(\ref{eq:markov_automaton}), and let $\pi$ denote the steady state distribution of the chain. 
Perfect Sampling Algorithm (PSA) gives a sample from the steady state distribution in finite time, using a coupling from the past construction.

\begin{theorem}[\cite{propp96exact}] 
\label{thm:psa}
Let $(a_{-t})_{t\in\N}=(a_0,a_{-1},\dots,a_{-t},\dots)$ be a sequence of \iid events with distribution $\init$ on $\evtSet$. There exists  $\ell \in \N$ such that 
$\lim_{t\to \infty} \big| \ph(\SS, a_{-t+1\to 0})
\big| =\ell \; \textrm{almost surely}.$
The grand coupling generated by $\M$ and $(a_{-t})_{t\in\N}$ is coalescing if $\ell = 1$. 
In that case, let 
$$\tau \bydef \inf \left\{ t \; : \; \big| \ph(\SS, a_{-t+1\to0}) \big| = 1 \right\}$$
be the coupling time of the chain. Then $\Eb(\tau) < \infty $ and 
$\ph(\SS, a_{-\tau+1 \rightarrow 0})$ is steady state distributed. 
\end{theorem}

The main drawback of PSA  is the fact
that one needs to simulate one Markov chain starting from each state in $\Space$, 
which is too large for most applications. 
Several approaches have been used to overcome this problem.
The main one for a 
partially ordered state space $(\Space, \preceq)$ and monotone events
was already given in \cite{propp96exact}.  
\begin{definition}
An event $a \in \Evt$ is said to be monotone if, for all $x,y\in\cX$, $x \preceq y \;
\Rightarrow \ph(x,a) \preceq \ph(y,a)$. 
\end{definition}
If all events are monotone, then one can consider only the trajectories issued from maximal and
minimal initial states (\cite{propp96exact}). 
In the case of general non-monotone chains, it is possible to use a bounding chain method, 
introduced in \cite{KendallMoeller-2000}. 
EPSA (Envelope Perfect Sampling Algorithm, \cite{BGV08,ashes}) constructs bounding chains in the case when the state space is equipped with a 
lattice order relation. We give next a short overview of EPSA. 

\subsection{Bounding Interval Chains}

Let $(\Space, \peq)$ be a 
lattice. For $m, M \in \Space$, denote by 
$[m, M] \bydef \{x \in \Space \; : \; m \peq x  \peq  M\}$
the \textit{(lattice) interval} between the endpoints $m$ and $M$.
Let $\J$ be the set of all nonempty lattice intervals:
$\J = \left \{[m, M] \; : \; m, M \in \Space, \; m \peq  M \right \}$.   
  Given a grand coupling $\{ (X_t(x))_{t \ge 0} \mid x \in \Space \}$, 
  a {\em bounding interval chain}  is
  any Markov chain of nonempty intervals $([m_t,M_t])_{t \ge 0}$
  such that: for all $x$ in $\Space$ and all 
  $t \ge 0$, $X_t(x) \in [m_t,M_t]$.
In particular we notice that when $m_t=M_t$, the grand coupling 
has necessarily coalesced. 

An \textit{envelope transition} function $\Ph : \J \times \Evt \to \J$
is defined by: for all $[m, M] \in \J$ and $a \in \Evt$,
\begin{eqnarray} \label{def:env_trans_fct}
       \Ph([m, M], a) 
                      \bydef \left [ \inf_{m \peq x\peq M } \ph(x, a), \;   
          \sup_{m\peq x\peq M} \ph(x, a) \right ].
\end{eqnarray}

As with $\ph$, the transition function $\Ph$ can be extended to finite words of events.
For any $t \in \N$, $F :{\J \times  \Evt^t} \to \J$ is defined  by:
$\Ph \left([m,M], \epsilon \right) \bydef [m,M]$ and $\Ph \left([m,M], a_{1\to t} \right) \bydef \Ph \left(\Ph \left([m,M], a_{1\to t-1} \right), a_{t} \right)$ for $t\geq 1$.

Let $\bot \bydef \inf \Space$ (resp . $\top \bydef  \sup
\Space$) be the bottom  
(resp.  top) element of $\Space$.
The process $[m_t,M_t] \bydef \Ph([\bot,\top], a_{1 \to t})$ is a Markov chain
over the state space   $\Space \times \Space$, called the {\em envelope chain},
and is a bounding interval chain of the grand coupling 
$\{(\ph(x,a_{1 \rightarrow t}))_{t\geq 0} \mid x \in \Space\}$.

The envelope process can be used to detect the coalescence of the grand coupling.
The following result was shown in \cite{BGV08}:
\begin{theorem}
\label{th:env}
Let $(a_{-t})_{t\in\N}$ be a sequence of \iid events with distribution $\init$ on $\evtSet$.
Assume that the envelope chain $\Ph([\bot,\top], a_{-t+1 \to 0})$ hits the set of single point intervals
$\mathcal{P} = \left \{[x,x] \; : \; x \in \Space \right \}$ a.s. 
in finite time. Let
$
{\tau_e} \bydef \min \left\{ t \; : \; \Ph([\bot,\top], a_{-t+1\to0}) \in \mathcal{P}\right\}, 
$
then $\tau_e$  is a  backward coupling time of the envelope
chain. The state defined by  $\Ph([\bot,\top], a_{-\tau_e+1\to 0})$
has the  steady state distribution of DTMC $(X_t)_{t\in \integer}$.
\end{theorem}

 \begin{algorithm}[H]
\KwData{
\Iid events  $\left(a_{-t}\right)_{t\in \integer} \in \Evt^\N$
}
\KwResult{A state $x^* \in \Space$ generated according to the stationary
 distribution of the Markov chain}
  \Begin{
$t:=1$\;
  \Repeat{$m = M$}
{ $m :=   \bot$; $M :=  \top$\;
 \For{$i=t-1$ {\em  \textbf{downto}} $0$}
{ $[m, M] :=\Ph \left( [m, M],a_{-i}\right)$ \;}
$t:=2t$\;}
$x^* := m$\;
 \Return{$x^*$}\;}
\caption{Envelope Perfect Sampling Algorithm (EPSA)}
\label{al:EPSA}
\end{algorithm}

Envelope Perfect Sampling Algorithm (EPSA)
is given in Algorithm~\ref{al:EPSA}. 
The reason to double $t$ at each iteration of the algorithm is that
we need to compute $\Ph([\bot,\top], a_{-t+1 \rightarrow 0})$ in each loop,
which corresponds to $t$ iterations of $\Ph$.
While increasing $t$ by $1$ would lead to a quadratic cost
in $\tau_e$, doubling it keeps the complexity linear. This was already observed in \cite{propp96exact}, for the monotone case.

The construction of the envelope chain depends on the
discrete event representation of the Markov chain that is not
unique. Different event representations lead to
different envelope chains with different coupling properties (one may coalesce
almost surely and the other not, or if they both coalesce their coupling times
may be different). The complexity of the envelope transition function may also differ
depending on the representation.


\begin{remark}
\label{rem:epsa}
When the assumptions of Theorem \ref{th:env} do not hold (i.e. the envelope chain does not couple), EPSA never stops. 
However in that case, the variants of EPSA can still provide perfect samples or performance bounds for increasing cost functions:
\begin{itemize}
 \item EPSA with splitting, proposed by \cite{BGV08}, can still generate perfect samples in finite time, under milder assumptions. 
The splitting algorithm is hybrid: it first runs EPSA when the envelopes are too far apart and switches to
the usual PSA algorithm as soon as the number of states inside the envelopes becomes manageable. 
 \item If we stop the algorithm at any time (for instance after some maximal number of iterations is reached), we will obtain an interval estimate for the stationary distribution.
Indeed, let $X \sim \pi$ be a random variable distributed according to the stationary distribution $\pi$. 
For any $s \geq 0$, $[\bot_s,\top_s] \bydef \Ph([\bot,\top], a_{-s+1\to 0})$ satisfies:
\begin{equation}
\label{eq:st}
\bot_s \st X \st \top_s, 
\end{equation}
where $\st$ denotes the usual strong stochastic order of random variables. For two random variables $X$ and $Y$ with values in $(\Space, \preceq)$, 
$X \st Y$ if $P(X \in U) \leq P(Y \in U)$ for all increasing sets $U \subset \Space$ (see \cite{Stoyan} 
for further material on stochastic orderings). 

This interval estimate can be used to obtain performance bounds, as explained in the following subsection. 
\end{itemize}

\end{remark}

\subsection{Performance Evaluation Bounds}
\label{ss:bounds}

Assume now $c : \Space \rightarrow \R$ is some increasing cost function defined on states of the Markov chain (e.g. the total number of items in replenishment, or the probability that a new demand cannot be fulfilled). 
Then (\ref{eq:st}) implies: 
$$
\Eb [c(\bot_s)] \leq \Eb[c(X)] \leq \Eb[c(\top_s)], \; s \geq 0,
$$
so we can use EPSA to obtain bounds of the steady-state cost. In performance evaluation of a given ATO system, or comparison between two different system designs, we are often interested 
in some specific performance guarantees (e.g. probability that a new demand cannot be fulfilled must be below a certain level). In that case,   
we are only interested in bounds for a given cost function and not the exact samples from the stationary distribution. Furthermore, the difference between the lower and upper bound provides also the error estimate for these bounds.
 
\medskip

In Sections \ref{sec:POS} and \ref{sec:TOS} we apply these existing methods (PSA and EPSA) to ATO systems with individual replenishments (PSA for the POS case in Section \ref{sec:POS}
and EPSA for the TOS case in Section~\ref{sec:TOS}).

\section{ATO-POS with individual state-dependent replenishments}
\label{sec:POS}

\subsection{Model description}
\label{subs:POS_state_dependence}

We consider the ATO-POS system with individual replenishments of items.  
There are $I$ different item types and let $\I = \{1, \ldots, I\}$. We assume finite stock capacities
and denote by $C_i$ the total amount of items of type $i \in \I$. 
Customers arrive in the system according to a Poisson process of rate $\lambda$. 
Each customer asks for a subset of items and the probability to ask for subset $A$ is denoted by $p_A$, \ie 
the demands for each subset $A$ follow a Poisson process of rate $\lambda_A = p_A \lambda$. 
If some demanded items are not available, then the customer takes the available items (POS case) and the demand for the items that are not available is lost. 
As often considered in the ATO literature, we assume that the number of different subsets customers can ask for is small (for instance, $|\{A\subset \I:\lambda_A\neq 0\}|$ is linear with respect to the number $I$ of item types).
Each item of type $i$ is replenished after an exponential time, with a rate that depends on the current amount of items $i$ in replenishment. We assume that the replenishments of different item types are mutually independent, 
and independent from the demands.

This system can be modeled as a network of $I$ queues with joint arrivals and independent services: arrivals to queues represent demands for different subsets of items and services in a queue model replenishments of items.
Denote by $C = (C_1, \ldots, C_{I})$ the vector of queue capacities. 

The total number of items in each queue (i.e. in replenishment) is given by a vector $x = (x_1, \ldots, x_I)$, where $x_i$ is the number of items of type $i$. 
The state space of the system is: $\cX=\{0,...,C_1\}\times \{0,...,C_2\} \times \cdots \times \{0,...,C_I\}.$ 
In the following, for $i\in \I$, we denote by $e_i$ the state with all the components equal to $0$,  
except component $i$ that is equal to~$1$. 

We have two different types of transitions. For each $x\in\cX$, $A \subset \I$ and for each $i\in \I$:
\begin{itemize}
 \item There is a demand for subset $A$, with rate $\lambda_A$. The new state is:
$
x+\sum_{j\in A}\mathds{1}_{\{x_j<C_j\}}e_j.
$

 \item If $x_i > 0$,  there is a service in queue $i$, with rate $\mu_i(x_i)$ that depends on the current number $x_i$ of items of type $i$. The new state is $x-e_i$.
\end{itemize}

By a standard uniformization procedure, we can transform this continuous time Markov chain to a discrete time Markov chain. 
Let $\beta_i:= \max_{1\leq x_i\leq C_i} \mu_i(x_i)$ be the maximal service rate for queue $i$.  
Then the outgoing rate in each state is upper-bounded by
$
\Lambda:=\lambda + \sum_{i \in \I} \beta_i.
$
We take the uniformization constant equal to $\Lambda$. 

\subsubsection*{Event representation.}
We now explain a discrete event representation of our (uniformized) Markov chain. 
In order to allow the construction of a family of Markov chains on the same probability space and driven by the same sequence of events, the set of possible events and the event rates cannot depend explicitely on the state. The following events allow such a construction.

\begin{itemize}
 \item \textit{Arrivals.}
For any $A \subset \I$, $A \neq \emptyset$, let 
$d_{A}$ be the event of probability $\lambda_{A}/\Lambda$ that corresponds to a ``joint arrival to queues in $A$''. 
 \item \textit{Services.} As the service rate in each queue does not depend on the state of other queues, we can consider the queues separately. 
 For queue $i$, the service rate is given by a function $\mu_i(x_i)$. We start by reordering the set $\{0, \ldots, C_i\}$ of possible values of $x_i$ in increasing order of $\mu_i$, 
 and denote this permutation by $\ell_i=(\ell_i^{(0)},...,\ell_i^{(C_i)})$. We have:
$$
 0=\mu_i(\ell_i^{(0)}) \leq \mu_i(\ell_i^{(1)}) \leq ... \leq \mu_i(\ell_i^{(C_i)})=\beta_i.
$$
For each $1\leq j\leq C_i$, let $s_i^{(j)}$ be the event of probability $\left[\mu_i(\ell_i^{(j)})-\mu_i(\ell_i^{(j-1)})\right]/\Lambda$ that corresponds to 
a ``service in queue $i$ for states $x=(x_1,...,x_I)$ such that $x_i\in \{\ell_i^{(j)},\ell_i^{(j+1)},...,\ell_i^{(C_i)}\}$''.
\end{itemize}

\begin{example}
We assume that $\mu_i(x_i)=\mu_i \cdot x_i$ for some $i\in \I$, and consider the events for services in queue $i$. 
For $1\leq j\leq C_i$, we define the event $s^{(j)}_i$ of probability $\mu_i/\Lambda$, as a service in the $i$-th queue for all the states $x$ such that $x_i\geq j$. 
In particular, states $x$ such that $x_i=C_i$ (that have the highest service rate in queue $i$) are served in each of these events, while, for states $x$ such that $x_i=1$, 
the number of items in the $i$-th queue  decreases only when event $s^{(1)}_i$ occurs. 
For $I=2$, $i=1$ and $C_1=4$, the transition function for events $s^{(3)}_1$ and $s^{(4)}_1$ is given in Figure \ref{fig:Services}.

\begin{figure}
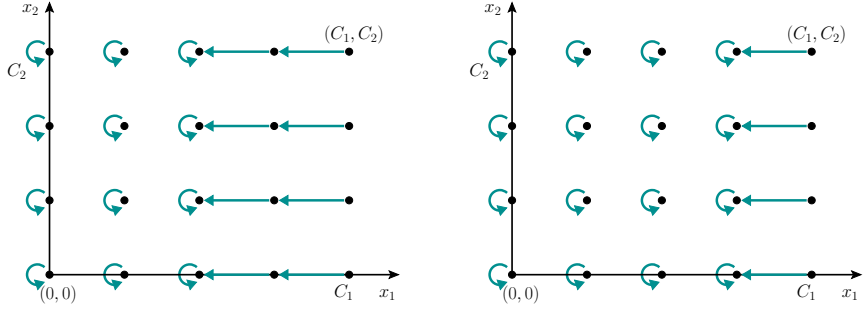

\begin{center}
\includegraphics[width=0.4\textwidth]{Services3.fig}
\includegraphics[width=0.4\textwidth]{Services4.fig} 
\caption{\label{fig:Services} Case $\mu_1(x_1)=\mu_1 \cdot x_1$. On the left: Transition function for service $s^{(3)}_1$.
On the right: Transition function for service $s^{(4)}_1$.} 
\end{center}
\end{figure}
\end{example}


The set of events in the system is:  $\cE = \{d_A, A \subset \I, A \neq \emptyset\} \cup \{ s_i^{(j)}, i\in \I, 1\leq j\leq C_i\}$
and the probability distribution $\nu$ on $\cE$ is given by: 
$\nu(d_{A}) = \lambda_{A}/\Lambda$, $\emptyset \neq A \subset \I$, $A \neq \emptyset$, and $\nu(s_i^{(j)}) = \left[\mu_i(\ell_i^{(j)})-\mu_i(\ell_i^{(j-1)})\right]/\Lambda$, 
$i\in \I, 1\leq j\leq C_i$. Note that some events may have probability $0$; we can ignore these events  (\ie restrict $\cE$ to the support of 
$\nu$).

\subsubsection*{Transition function.}
The transition function $\fp : \mathcal{X} \times \cE \rightarrow \mathcal{X}$ of the ATO-POS system is defined as follows. 
For $x\in\mathcal{X}$, $A\subset \I$, $A\neq\emptyset$, $i\in \I$ and $1\leq j\leq C_i$:
\begin{eqnarray*}
\left\{ \begin{array}{rcl} 
        \fp(x,d_A) &=& x+\displaystyle\sum_{k\in A}\mathds{1}_{\{x_k<C_k\}}e_k , \\
	\fp(x,s_i^{(j)}) &=& x-\displaystyle\sum_{k=j}^{C_i}\mathds{1}_{\{x_i=\ell_i^{(k)}\}}e_{i} \;\: = \;\: x-\mathds{1}_{\left\{\mu_i(x_i)\geq \mu_i\left(\ell_i^{(j)}\right)\right\}}e_{i} .
        \end{array}
\right.
\end{eqnarray*}

\subsubsection*{Monotonicity.}
We consider the natural product order on $\mathcal{X}$, \ie $x=(x_1,...,x_I) \leq (y_1,...,y_I)=y$ if $x_i\leq y_i$ for all $1\leq i\leq I$. 
By using the fact that a service $s_i^{(j)}$ can only modify the $i$-th coordinate, it easily follows that:
\begin{proposition}
\label{prop:ATOpos_monotone}
The transition function $\fp$ of the ATO-POS system is monotone 
under the natural product order on $\mathcal{X}$, \ie for all $a \in \cE$,
$x \leq y \; \Longrightarrow \fp(x,a) \leq \fp(y,a).$
\end{proposition}

\smallskip

Since the system is monotone, it is enough to consider only two trajectories (starting from the upper and lower states). 
The complexity of the Perfect Sampling Algorithm depends on the coupling time of these two trajectories. 

\subsection{Bound for the coupling time}
\label{subsubs:POS_slightly_general}

Let $(a_{-t})_{t\in\N}=(a_0,a_{-1},\dots,a_{-t},\dots)$ be a sequence of \iid events with distribution $\init$ on $\evtSet$.
Let $\tau^P:=\min \left\{ t \; : \; \big|\fp(\cX, a_{-t+1\to0})\big| =1\right\}$ be the coupling time of the whole chain $X$ (defined on $\cX$). For all $1\leq i\leq I$, we define the projection on $i$-th component as: 
\[\phi_i:\left\{\begin{array}{ccc} \cX & \longrightarrow & [0,C_i] \\ x=(x_1,\dots,x_I) & \longmapsto & x_i \\ \end{array} \right.,\]
and set $\tau^{P,i}:=\min \left\{ t \; : \; \big|\phi_i\left(\fp(\cX, a_{-t+1\to0})\right)\big| =1\right\}$, the ``coupling time'' of the chain on the $i$-th component. 

We first prove the following property:
\begin{eqnarray}\label{eq:property}
\Big( \: \big|\phi_i\left(\fp(\cX, a_{-s+1\to0})\right)\big| =1 \textrm{ and } t\geq s \: \Big) \; &\Longrightarrow& \;\big|\phi_i\left(\fp(\cX, a_{-t+1\to0})\right)\big| =1.
\end{eqnarray}
Informally, we say that, as soon as a component couples, it \textit{stays coupled}. This property implies that $\tau^P=\max_{1\leq i\leq I} \tau^{P,i}$.

\begin{definition} 
\label{def:cc}
A Markov chain that satisfies property \eqref{eq:property} for any sequence of events $(a_{-t})_{t\in\N}$ is said to be \textit{componentwise coupling}. In that case, the whole coupling time of the chain can be expressed in terms of the coupling times on each component.
\end{definition}

The following lemma shows that the chain $X$ is componentwise coupling.
\begin{lemma} \label{lem:property}
Let $i\in\I$. The projection $\phi_i(X)$ of the chain on the $i$-th component is a Markov chain on the state space $[0,C_i]$. Let $f_{P,i}$ be its transition function. Then we have that: 
\[\tau^{P,i}=\min \left\{ t \; : \; \big|f_{P,i}([0,C_i], a_{-t+1\to0})\big| =1\right\},\]
\ie $\tau^{P,i}$ is the coupling time of the Markov chain $\phi_i(X)$ (it does not depend on the value of the chain $X$ on other components).
\end{lemma}

\proof 
 We compute $\phi_i\circ\fp$. For $x\in\mathcal{X}$, $A\subset \I$, $A\neq\emptyset$, $i\in \I$ and $1\leq j\leq C_i$:
\begin{eqnarray*}
\left\{ \begin{array}{rcl} 
        \phi_i\left(\fp(x,d_A)\right) &=& x_i+\mathds{1}_{\{i\in A \:\&\: x_i<C_i\}} , \\
	\phi_i\left(\fp(x,s_i^{(j)})\right) &=& x_i-\mathds{1}_{\left\{\mu_i(x_i)\geq \mu_i\left(\ell_i^{(j)}\right)\right\}} .
        \end{array}
\right.
\end{eqnarray*}
Hence we have that $\phi_i\circ\fp$ only depends on $x_i$. For any event $a\in\cE$, we can set $f_{P,i}(x_i,a) = \phi_i\left(\fp(x,a)\right)$. This ends the proof.
\endproof

As a direct application of Lemma \ref{lem:property}, we have:
\begin{eqnarray}
\label{eq:sum_tau}
 \Eb[\tau^P] = \Eb\left[\max_{1\leq i\leq I} \tau^{P,i}\right] \leq \sum_{i=1}^I \Eb[\tau^{P,i}].
\end{eqnarray}
 In order to bound $\Eb[\tau^P]$, we will show bounds on $\Eb[\tau^{P,i}]$.
We use the following result that is often used as a part of the proof of Foster's stability criterion (see for instance \cite[proof of Theorem~1.1]{Bremaud}). 
\begin{lemma}
\label{Foster}
 Let the transition matrix $P$ on the finite state space $S$ be irreducible and suppose that there exists a function $h:S\rightarrow \mathbb{R}_+$ such that
\begin {eqnarray}
\label{eq:Foster_cond}
 \sum_{z\in S} P(y,z) h(z) \leq h(y)-\epsilon  \textrm{ for all } y\notin U,
\end {eqnarray}
for some subset $U\subset S$. Let $\tau_U$ be the hitting time of $U$ and $\Eb_y$ denote the expectation, knowing that the chain starts in $y$. Then, for all $y\notin U$,
\begin{eqnarray}
\label{eq:Foster_ccl}
 \Eb_y[\tau_U] \leq \frac{h(y)}{\epsilon}.
\end{eqnarray}
\end{lemma}

We next give a bound for $\Eb(\tau^{P,i})$ for the two following cases: either service or arrival rate is high.
\begin{lemma}
\label{lem:tau}
 Let $i\in\I$. Let $\lambda_i:=\sum_{A\ni i} \lambda_A$ be the total arrival rate in queue $i$. Set $\alpha_i :=\min_{1\leq x_i\leq C_i} \mu_i(x_i)$ and $\eta_i :=\max_{1\leq x_i < C_i} \mu_i(x_i)$. Then: 
\begin{itemize}
 \item If $\lambda_i < \alpha_i$, then $\Eb[\tau^{P,i}] \leq\frac{\Lambda C_i}{\alpha_i-\lambda_i}$.
 \item If $\lambda_i > \eta_i$, then $\Eb[\tau^{P,i}] \leq\frac{\Lambda C_i}{\lambda_i-\eta_i}$.
\end{itemize}
\end{lemma}

\proof 
 Let $S=\{0,...,C_i\}$ and $P$ be the transition matrix of the Markov chain $Y=\phi_i(X)$. When the chain $Y$ is in state $y$, $0 \leq y \leq C_i$: $Y$ goes to $y+\mathds{1}\{y<C_i\}$ with probability $\lambda_i/\Lambda$, 
it goes to $y-1$ with probability $\mu_i(y)/\Lambda$, and stays at $y$ with probability $\left[\Lambda-\lambda_i-\mu_i(y)\right]/\Lambda$. 
 
In the first case ($\lambda_i <\alpha_i$), we use Lemma \ref{Foster} with $U=\{0\}$, $h(z)=z$ for all $z\in S$, $\epsilon=(\alpha_i-\lambda_i)/\Lambda$ and apply \eqref{eq:Foster_ccl} with $y=C_i$. This gives that $\Eb[\tau^{P,i}]\leq \Eb_{C_i}[\tau_{\{0\}}] \leq \frac{h(C_i)}{\epsilon}=\frac{\Lambda C_i}{\alpha_i-\lambda_i}$, where the first inequality comes from the fact that the chain $X$ (and thus its projection $\phi_i(X)$) is monotone (Proposition \ref{prop:ATOpos_monotone}): when the chain starting from $y=C_i$ reaches $0$, the chains starting from all other states are in $0$, and the system coupled. 
In the second case ($\lambda_i>\eta_i$), we use Lemma \ref{Foster} with $U=\{C_i\}$, $h(z)=C_i-z$ for all $z\in S$, $\epsilon=(\lambda_i-\eta_i)/\Lambda$ and apply \eqref{eq:Foster_ccl} with $y=0$.
\endproof

This lemma and equation \eqref{eq:sum_tau} lead to the following proposition:
\begin{proposition}
\label{prop:POS_dim}
 Let $\I_0=\{i:1\leq i\leq I, \lambda_i < \alpha_i \}$, $\I_C=\{i:1\leq i\leq I, \lambda_i > \eta_i\}$ and assume $\I=\I_0 \cup \I_C$. Then we have:
$$
 \Eb[\tau^P] \leq \Lambda\left(\sum_{i\in \I_0} \frac{C_i}{\alpha_i-\lambda_i} + \sum_{i\in \I_C} \frac{C_i}{\lambda_i-\eta_i}\right).
$$
\end{proposition}

We discuss two important cases:
\begin{itemize}
 \item {\em Single server case.} Assume $\mu_i(x_i)=\mu_i$ for all $i\in\I$, $x\in\cX$. Then $\alpha_i=\eta_i=\mu_i$, and $\Lambda=\lambda+\sum_i\mu_i$. 
Hence the mean coupling time of the chain $X$ is in $O\left(|C|\right)$, with $|C|=\sum_i C_i$. In addition we can notice that the hypothesis $\I = \I_0 \cup \I_C$ is necessary to have a linear bound in $|C|$: 
$\I \neq \I_0 \cup \I_C$ implies $\lambda_i=\mu_i$ for some $i\in\I$, and $\Eb[\tau^{P,i}]$ is quadratic in $C_i$ in that case (see \cite{Dopper}).
\item {\em Infinite server case.} Assume $\mu_i(x_i)=\mu_i \cdot x_i$ for all $i\in\I$, $x\in\cX$. 
Then $\alpha_i=\mu_i$, $\eta_i=\mu_i \cdot (C_i-1)$ and $\Lambda=\lambda+\sum_i\mu_i\cdot C_i$. 
If $\I=\I_0$, the mean coupling time of the chain is in $O\left(|C|^2\right)$. Note that this bound is larger than in the single server case as a result of the time-discretization (uniformization) of our chain.
\end{itemize}

In Appendix \ref{sec:AppB}, we consider a slightly more general model: the service rates $\mu_i$ can depend on other components, but the service events remain monotone. 
Under a high service rate assumption, we give a bound on the mean hitting time to zero (that provides an upper bound for the coupling time). 
This result will be also used to bound the running time of the AEPSA algorithm for ATO-POS with joint returns, studied in Section \ref{sec:Aggreg_POS}.

\section{ATO-TOS with individual state-dependent replenishments}
\label{sec:TOS}

\subsection{Model description}

We consider the ATO-TOS model with individual replenishments of items. The difference with the ATO-POS 
model is in the way the demands are handled in the out-of-stock situation:
If some demanded items are not available, then the whole demand is lost. 

As before, we model the system by a queueing system, with state space $\cX$. 
We consider the uniformized Markov chain with uniformization constant $\Lambda=\lambda + \sum_{i \in \I} \beta_i$, with $\beta_i= \max_{1\leq x_i\leq C_i} \mu_i(x_i)$ as before. We consider the same set $\cE$ of events and the same probability distribution $\nu$ on $\cE$ as in the POS case (defined in Section~\ref{subs:POS_state_dependence}). 
The transition function for services is also the same as in the POS case.
For an arrival $d_A$, $A\subset \I$, $A\neq \emptyset$, the response of POS and TOS systems is different only for states $x$ that belong to the boundary of the state space (\ie such that there exists $i\in A$ with $x_i=C_i$).

\subsubsection*{Transition function.}
The transition function $\ft$ for the ATO-TOS system is defined as follows. Let $A\subset \I$, $A\neq\emptyset$, $i\in \I$, $1\leq j\leq C_i$, and $x\in\mathcal{X}$:
\begin{eqnarray*}
\left\{ \begin{array}{rcl} 
        \ft(x,d_A) &=& x+\displaystyle \left(\prod_{i\in A}\mathds{1}_{\{x_i<C_i\}}\right) e_{A} , \\
	\ft(x,s_i^{(j)}) &=& \fp(x,s_i^{(j)}),
        \end{array}
\right.
\end{eqnarray*}
where $e_A=\sum_{i\in A}e_i$. In other words, when an arrival $d_A$ occurs, we add 1 to each component of $x$ in $A$, if \textit{all} its components $x_i$, $i\in A$, satisfy $x_i<C_i$. 

In the ATO-TOS system, whether or not an arrival in queue $i$ is accepted depends on the whole state of the system, which makes the system more difficult to study. 
Indeed, contrary to ATO-POS, arrivals are not necessarily monotone in the ATO-TOS system.

\subsection{Envelopes} \label{subs:envelopeTOS}

\subsubsection*{(Non-)monotonicity.}

As before, we can consider the product order on $\cX$.
Since services are the same as in the ATO-POS model, they are monotone (Proposition \ref{prop:ATOpos_monotone}). For the same reason, the arrivals of \textit{only one} item are also monotone.
Unfortunately, as soon as $|A|\geq 2$, the event $d_A$ is not monotone for the product order on $\cX$. For instance, let $I=2$, $x=(C_1-1, C_2-1)$ and $y=(C_1-1,C_2)$. Then $x\leq y$, yet $\ft(x,d_{\{1,2\}})=(C_1,C_2) \geq y=\ft(y,d_{\{1,2\}})$. 

One could try to find another partial order on $\cX$ for which the ATO-TOS model would be monotone. 
If $I=2$, it is easy to check that this is true for the following partial order: $x=(x_1,x_2)\preceq (y_1,y_2)=y$ if $x_1\geq y_1$ and $x_2\leq y_2$.
However, for $I\geq 3$ we show much stronger statement (the proof is given in Appendix \ref{sec:AppI3}):
\begin{proposition}
\label{prop:I3}
 Let us consider the ATO-TOS system with $I\geq 3$. If $\ft$ is monotone for a partial order $\preceq$, then $\preceq$ is the trivial order, \ie: $x\preceq y \Leftrightarrow x=y$.
\end{proposition}

We now show how to appy EPSA (Algorithm \ref{al:EPSA}) to the ATO-TOS model. 

\subsubsection*{Computation of envelopes.}

We consider the product order on $\cX$, and define the envelope transition function $F_T$ as in \eqref{def:env_trans_fct}, \ie for all $m$, $M\in\mathcal{X}$ such that $m\leq M$, and $a \in \Evt$, we set:
\begin{eqnarray*} 
       \Ph_T([m, M], a) 
                      \bydef \left [ \inf_{m \leq x\leq M } \ph_T(x, a), \;   
          \sup_{m\leq x\leq M} \ph_T(x, a) \right ].
\end{eqnarray*}
Services and arrivals of only one object are monotone for $\ft$, therefore envelopes follow easily: let $m$, $M\in\mathcal{X}$ such that $m\leq M$, and consider a monotone event $a$, then we have $ \Ft([m,M],a)=\left[ \ft(m,a),\ft(M,a) \right] $.

Thus we are left with the computation of envelopes for an arrival $d_A$ of several objects. As a consequence of the following proposition, we have that the computation of envelopes can be done in a linear time with respect to the number $I$ of queues.

\begin{proposition}
\label{prop:envTOS} 
Let $A\subset \I$ such that $|A|\geq 2$. Let $m$, $M\in\mathcal{X}$, $m\leq M$. We distinguish three cases:
\begin{description}
\item[Case I:] For all $i\in A$, $M_i<C_i$. Then
$ \Ft([m,M],d_A)=\left[m+e_A,M+e_A\right] $.

\item[Case II:] There exists $i_0\in A$ such that $m_{i_0}=C_{i_0}$. Then
$ \Ft([m,M],d_A)=\left[m,M\right]$.

\item[Case III:] Otherwise, \ie if there exists $i_0\in A$ such that $M_{i_0}=C_{i_0}$ and for all $i\in A$, $m_i<C_i$.
Let $[m',M'] := \Ft([m,M],d_A)$.
Then:
\begin{eqnarray*}
 m'_{i_0} &=& m_{i_0}+\prod_{k\in A\setminus\{i_0\}}\mathds{1}_{\{M_k<C_k\}}, \\
m'_i &=& m_i, \; i\not=i_0,
\end{eqnarray*}
and $M'_i = M_i+\mathds{1}_{\{i\in A ,\: M_i<C_i\}}$.
\end{description}
\end{proposition}

\proof 
Cases I and II are straightforward. In case I, for all $x\in \mathcal{X}$ such that $m\leq x\leq M$, we have $\ft(x,d_A)=x+e_A$. In case II, for all $x\in \mathcal{X}$ such that $m\leq x\leq M$, we have $\ft(x,d_A)=x$.

 {\em Case III.} 
 If $i\notin A$, then $m'_i=m_i$ and $M'_i=M_i$. Let $i\in A$. 

The upper envelope is simpler. Clearly, $M'_i\leq M_i+\mathds{1}_{\{M_i<C_i\}}$, 
so we only need to find $x$, $m\leq x\leq M$, such that $\left(\ft(x,d_A)\right)_i=M_i+\mathds{1}_{\{M_i<C_i\}}$. 
This is true for state $x$ such that $x_i=M_i$ and $x_j=m_j$ for all $j\neq i$ 
(by assumption, for all $j\in A$, $m_j<C_j$, so an arrival $d_A$ cannot be 'blocked' by $j\neq i$). 

Now we compute the lower envelope. Clearly, $m'_i \geq m_i$, for all $i$. 
For state $x$ such that $x_{i_0} = M_{i_0}$ and $x_i = m_i, \; i \not=i_0$, we get $\left(\ft(x,d_A)\right)_i=m_i, \; i\not=i_0,$ so 
$m'_i = m_i, \; i\not=i_0.$
For component $i_0$, we distinguish two cases:
\begin{enumerate}
 \item $\exists k\in A\setminus\{i_0\},\: M_k=C_k$. This case is similar as before, as we can take $x_k=M_k$ and $x_i = m_i, \; i\not=k$ ($x$ is 'blocked' by the $k$-th component, and $\ft(x,d_A)=x$). As $i_0 \not=k$, we get $m'_{i_0} = m_{i_0}.$ 
 
 \item $\forall k\in A\setminus\{i_0\},\: M_k<C_k$. 
 For all $x\in\mathcal{X}$ such that $x_{i_0}=m_{i_0}$, we have that $x_k\leq M_k<C_k$ for $ k\in A\setminus\{i_0\}$, so $\ft(x,d_A)=x+e_A$. 
 Hence $\left(\ft(x,d_A)\right)_{i_0}=m_{i_0}+1$. 
\end{enumerate}
\endproof

\subsection{Bound for the coupling time}

Let $\tau^T$ be the coupling time of the ATO-TOS system. Using a coupling between POS and TOS models and a bound on the mean hitting time to zero for the POS model give the following proposition:

\begin{proposition} \label{prop:bound_TOS}
 Assume $\delta=\min_{x\neq 0} \sum_i \mu_i(x)-\sum_i \lambda_i > 0$. Then we have:
$
 \Eb[\tau^T] \leq \frac{\Lambda}{\delta} |C|,
$
where $|C| = \sum_{i=1}^I C_i.$
\end{proposition}

Before proving this proposition, we remark that, as for Proposition \ref{prop:POS_dim}, the mean coupling time of the chain $X$ is in $O\left(|C|\right)$ for the single server case ($\mu_i(x_i)=\mu_i$ for all $i\in\I$), and in $O\left(|C|^2\right)$ when $\mu_i(x_i)=\mu_i \cdot x_i$ for all $i\in\I$ and $x\in\cX$.

In order to prove Proposition \ref{prop:bound_TOS}, we first need the following lemma:

\begin{lemma}
\label{lem:POS_bound_TOS}
 Let $a\in\cE$ and $x,y\in\cX$. Then $x\leq y \; \Rightarrow \; \ft(x,a)\leq \fp(y,a)$.
\end{lemma}

\proof[Proof of Lemma \ref{lem:POS_bound_TOS}.]
 Result for services follows from the fact that they are the same in both of the models, and from their monotonicity. In the case of an arrival $d_A$, $A\subset \I$, $A\neq \emptyset$, we have, for all $x\in\cX$: $\ft(x,d_A)\leq \fp(x,d_A)$. We can conclude using the monotonicity of arrivals for the ATO-POS model.
\endproof

\proof[Proof of Proposition \ref{prop:bound_TOS}.]
Let $\tau_0^P$ be the hitting time to zero for the ATO-POS model. We show that 
\begin{eqnarray} \label{eq:TOS_POS}
 \tau^T\leq \tau_0^P.
\end{eqnarray}
Then Proposition \ref{prop:POS_zero} (Appendix \ref{sec:AppB}) gives a bound on the mean of $\tau_0^P$ and concludes the proof. Now we prove \eqref{eq:TOS_POS}.
Assume we do coupling from the past, starting from state $C$, for both ATO-POS and ATO-TOS (with a coupling using the same events). When ATO-POS is in state zero, ATO-TOS is also in state zero by Lemma \ref{lem:POS_bound_TOS}, and hence it has coupled. 
So \eqref{eq:TOS_POS} follows.

Note that we cannot bound $\tau^T$ by the coupling time $\tau^P$ of ATO-POS (see Figure  \ref{fig:ATOtimes} and the following subsection).  
\endproof

\subsection{Comparison between POS and TOS models}

\begin{figure}
\begin{center}
\includegraphics[width=0.51\textwidth]{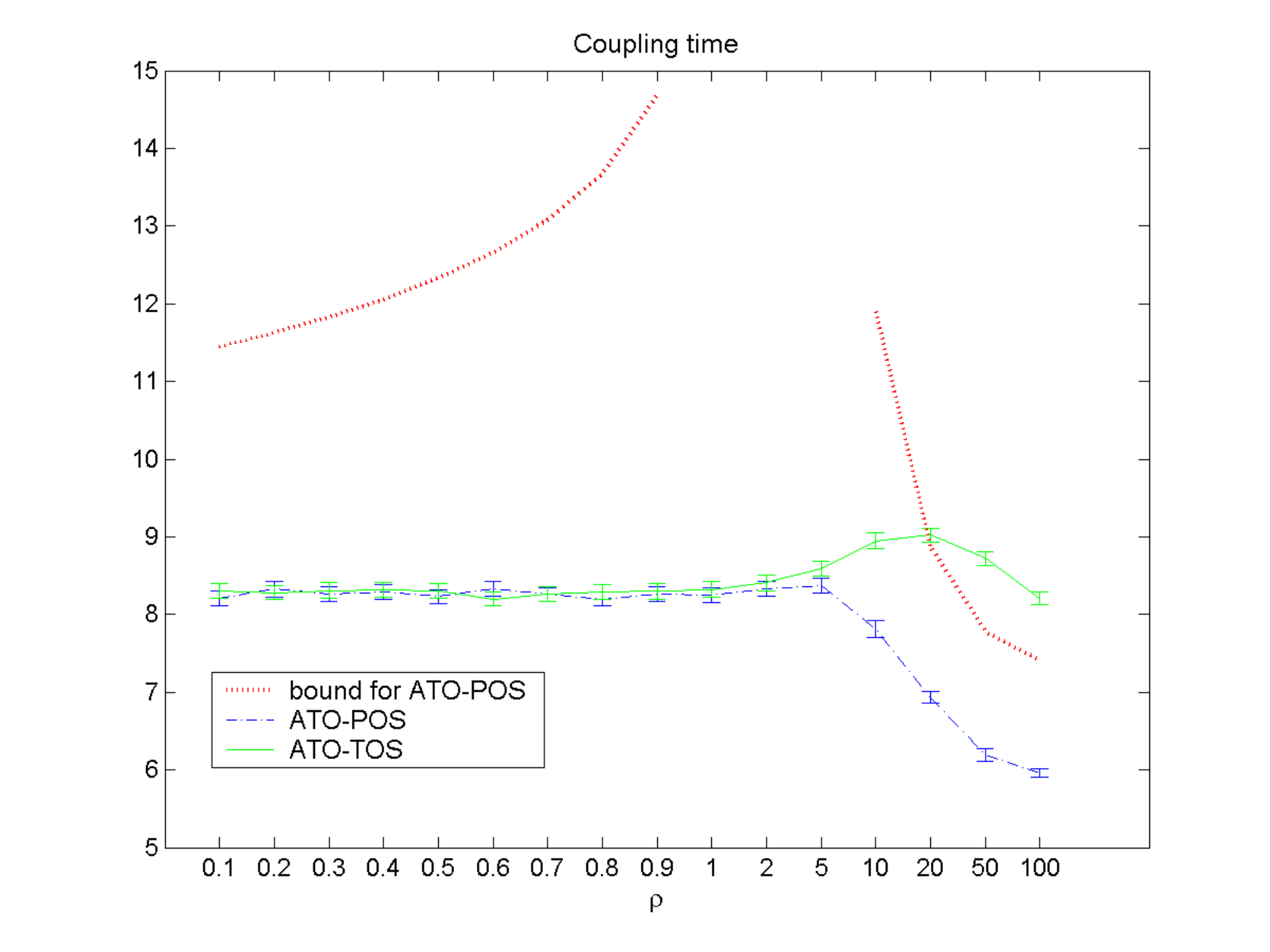}
\hspace{-2em}
\includegraphics[width=0.51\textwidth]{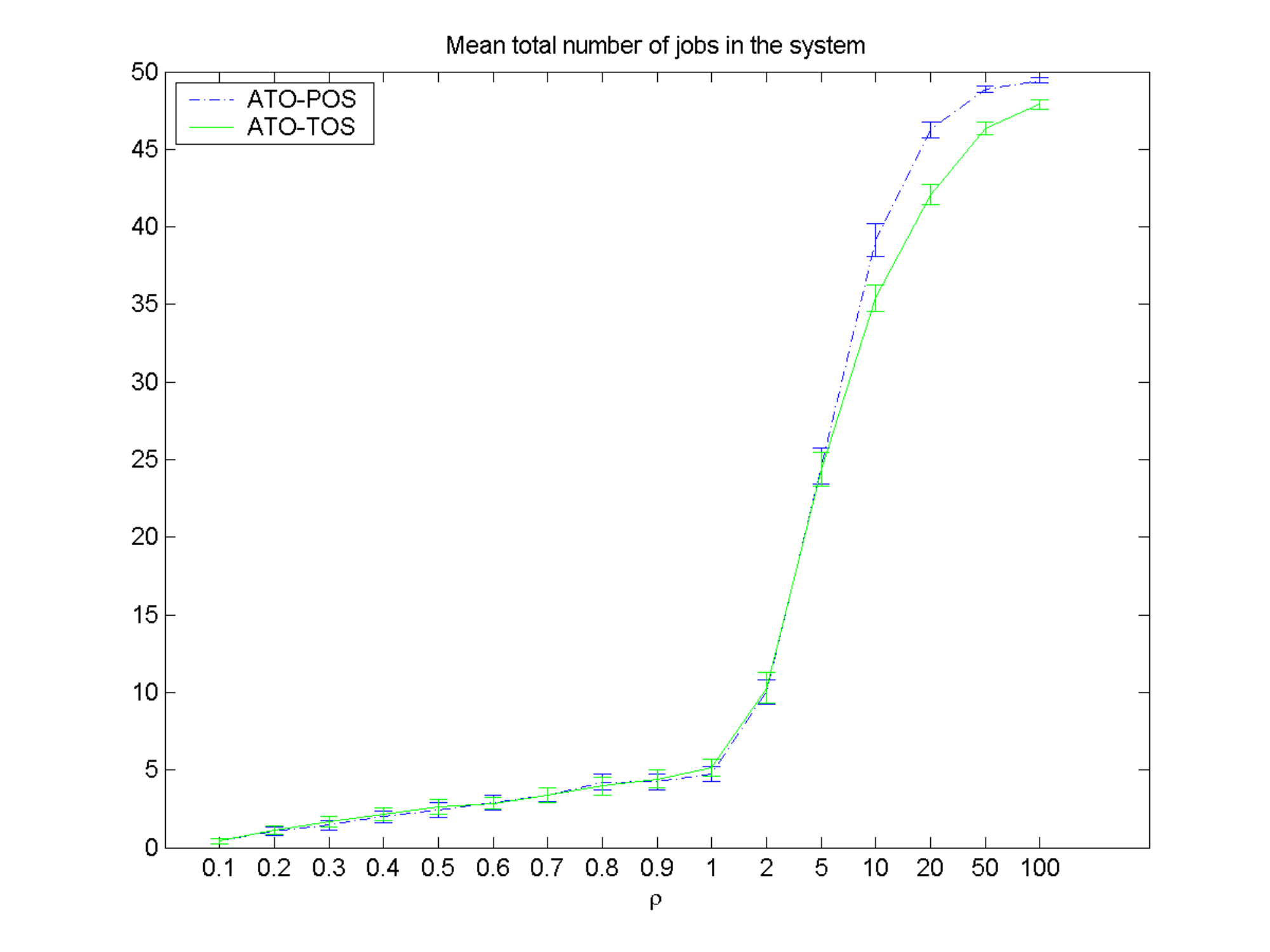}  
\caption{\label{fig:ATOtimes} On the left: Coupling times for ATO-POS and ATO-TOS, and the upper bound for ATO-POS (we display $\log_2(T)$ where $T$ is the mean coupling time), together with the $95\%$ confidence intervals.
On the right: Mean number of jobs in the system  for ATO-POS and ATO-TOS, together with the $95\%$ confidence intervals.} 
\end{center}
\end{figure}

In Figure \ref{fig:ATOtimes}, on the left, we give coupling times for ATO-POS and ATO-TOS models with individual services (Sections \ref{sec:POS} and \ref{sec:TOS}), for the following parameters: $I=5$, $C_i=10, \forall i$, $\lambda_A = \frac{1}{2^{|A|-1}}$, and $\mu_{i}(x_i) = \mu x_i$, with $\rho = \frac{\lambda_i}{\mu}$. The size of the sample is $n=100$. We can see that the upper bounds for ATO-POS (see Proposition \ref{prop:POS_dim}) are quite pessimistic for this example, in particular when $\lambda_i < \alpha_i$. This is mainly due to the fact that these bounds were obtained by considering the sum of the mean coupling times for each dimension, while the coupling of other components may occur faster knowing that one component has already coupled. In addition, when $\lambda_i >> \mu$, we can observe that the bound we obtained for the POS model is not a bound for the TOS model (even if the transition function of the POS model is  a bound for the one of the TOS model). This is an illustration of the fact 
that the chain in the TOS model is not componentwise coupling (see the definition above Lemma \ref{lem:property}, in \ref{subsubs:POS_slightly_general}). This explains why we use the hitting time to zero in the POS system, in order to bound the coupling time of the TOS chain.
 
In Figure \ref{fig:ATOtimes}, on the right, we give the mean values of the total number of jobs in the system (\ie $\sum_i x_i$) for ATO-POS and ATO-TOS, for the same parameter values. One can see that both models are very close when $\rho<1$, and the difference between the two models only become significant when $\lambda_i >> \mu$, which is expected as the difference between them is only for joint arrivals 
when the system reaches its capacity limits. 

In Section \ref{sec:Aggreg_POS}, we will consider ATO systems with joint services. We explain in \ref{subs:POS_joint_services} why the existing methods (PSA or EPSA) cannot be used to sample (with reasonable complexity) the steady state distribution. Hence we introduce a new method, that we explain in the following section.

\section{Aggregated envelopes}
\label{sec:method}

We introduce a new method of perfect sampling, that we apply in Section \ref{sec:Aggreg_POS} to an Assemble-To-Order system with Partial Order Service (ATO-POS) and joint returns of items, 
also called the service tools model by \cite{VH09}.
In this model, the customers demand (or borrow) subsets of items and return them together (see Section \ref{subs:POS_joint_services} for more details). 
Although we are usually interested only in the total number of available items of each type, this information is not sufficient to describe the evolution of the system: 
In order to get a Markov chain, we need to keep track of the way items leave the stock, as they will be returned together. 
The state space becomes rapidly intractable: 
Even the dimension of the state space is exponential with respect to the number $I$ of different item types. 
Thus even storing the vector representing the state of the system becomes challenging.

The idea of aggregated envelope method is to consider the projection of the state space on a more tractable space: In ATO-POS system with joint returns, we consider the projected space $\cX$ of vectors 
$x=(x_1,\dots,x_I)$, where $x_i$ is the total number of available type $i$ items. Considering space $\cX$ reduces exponentially the dimension of the state space.
However, vector $x$ does not contain all the information about the evolution of the system so we will need to construct a bounding chain that takes into account all possible evolutions.

More generally, we assume that our initial Markov chain is given by a tuple $\M = (\cS,\evtSet,\init,\ffp)$. 
We assume further that there is a projection function $\proj:\cS\rightarrow \cX$ such that $(\cX, \preceq)$ is a finite lattice. In practice, $|\cX|$ will be much smaller that $|\cS|$.
The state space $\cS$
is not necessarily a lattice (we do not assume any ordering relation on $\cS$). This is another important motivation for the aggregated envelope method: The state space of the ATO-POS system with joint returns 
is not a lattice for the product order.

In this section, we develop a method that samples an interval of $\cX$ containing the projection of a state distributed according to the stationary distribution $\pi$ of the original Markov chain. 
We will see that, under some conditions (see Section \ref{sec:POS_joint_returns}), it is even possible to sample a state in $\cS$, distributed according to the stationary distribution.

\subsection{Aggregation}

\begin{figure}
 \begin{center}
\scalebox{0.4}{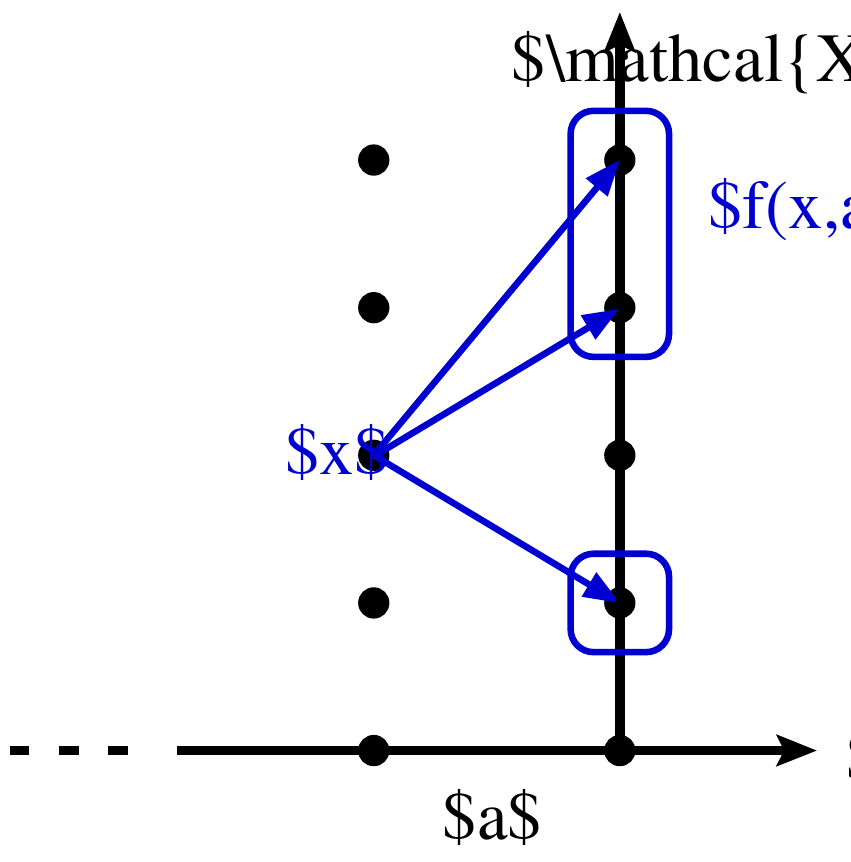}
\hspace{6em}
\scalebox{0.4}{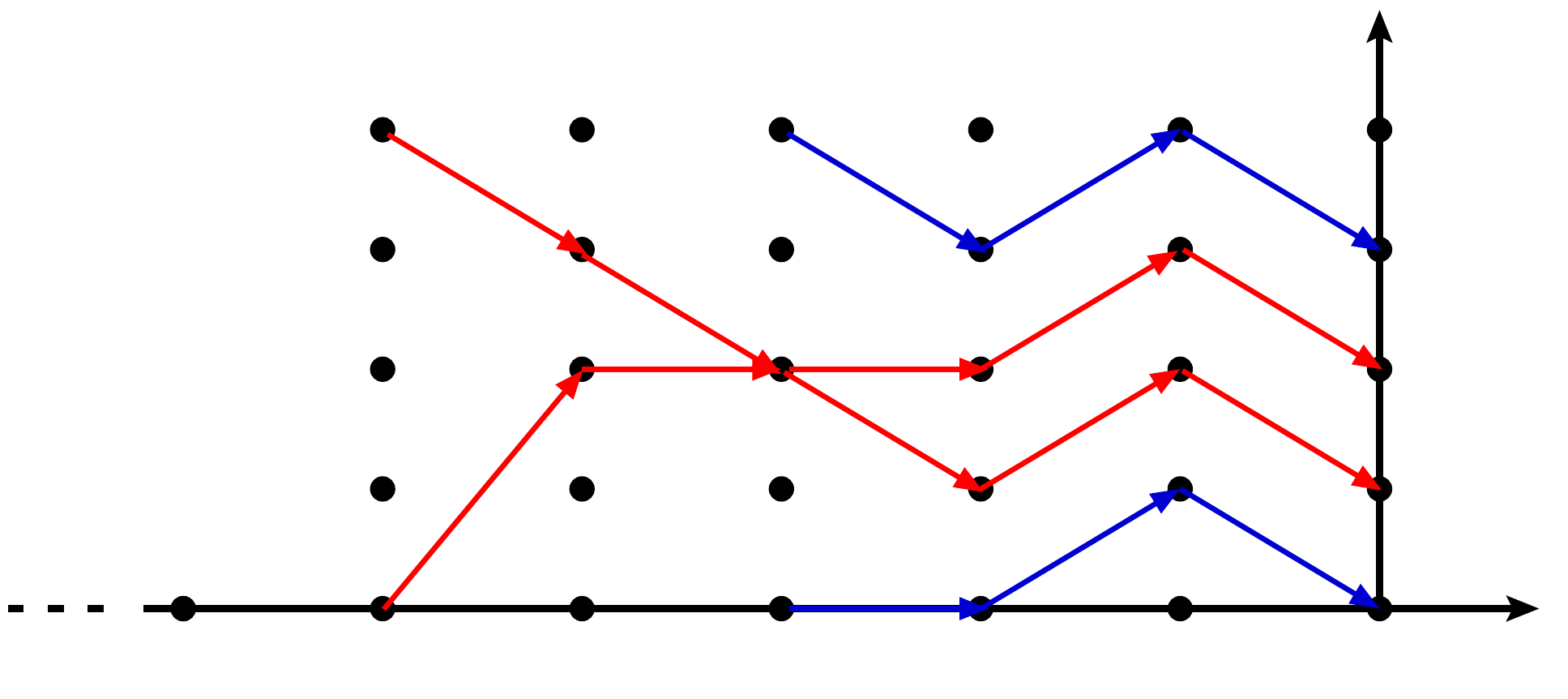}
\caption{\label{fig:def_f} On the right: The idea of aggregation (we consider the evolution from {\em all} the states $n \in \cS$ such that $\proj(n) = x$). On the left: Even if at some time $-t$ we have only one value for the projected process, this is not necessarily the case for times $-s$, $-t\leq -s \leq 0$, as  
this single projected value at time $-t$ can correspond to many different states in $\cS$.} 
\end{center}
\end{figure}

Our starting idea is to use the projected state space $\cX$ 
for simulations. Intuitively, the original Markov chain $X = (X_t)_{t\in \N}$ evolves in $\cS$, but we can only observe its projection $Y=(Y_t)_{t\in \N} = (\proj(X_t))_{t\in \N}$.
Assume that the original chain is in state $n_0 \in \cS$. The only information we have is $x = \proj(n_0) \in \cX$.  
When an event $a\in\cE$ occurs, we need to determine the next state in $\cX$. As $Y$ is not a Markov chain, 
we cannot determine the next state only from knowing $x$. Instead, we will consider the evolution from {\em all} the 
states $n \in \cS$ such that $\proj(n) = x$. More formally, for $x\in \cX$, we consider the 
following \textit{subset} $\Ux \subset\cS$: 
\begin{equation*}
\Ux=\{n\in\cS, \proj(n)=x\}=\proj^{-1}(\{x\}).
\end{equation*}
and define the following function $f : \cX \times \cE \rightarrow \PcX$, as illustrated in Figure~\ref{fig:def_f}:
\begin{equation}
f(x, a)=\proj\left(\ffp(\Ux,a)\right),
\end{equation}
where $\ffp$ is the transition function of the original Markov chain.

Let $\PcX$ denote the family of subsets of $\cX$, and 
$F : \PcX \times \cE \rightarrow \PcX$ be the transition function defined by:
\begin{equation}
\label{eq:agf}
F(U, a)=\cup_{x\in U} f(x,a).
\end{equation}
Since $F$ is a bounding chain for the projected process, we could use a coupling from the past scheme (for one trajectory starting from $U=\cX$) to provide an interval containing the projection of a state distributed according to the stationary distribution $\pi$: This is Algorithm~\ref{al:agg}. However, as explained below, this algorithm has several drawbacks and we will rather combine it with the method of envelopes developed in \cite{BGV08}.


\begin{algorithm}[ht]
\KwData{
\Iid events  $\left(a_{-t}\right)_{t\in \integer} \in \Evt^\N$.
}
\KwResult{Subset $U\subset\cX $ containing the projection of a state $n^*\in\cS$ distributed according to $\pi$.}
  \Begin{
$t=1$; $c=0$\;
  \Repeat{$c=1$}
{ $U =\cX$\; 
 \For{$i=t-1$ {\em  \textbf{downto}} $0$}
{ $U :=F\left(U,a_{-i}\right)$ \;
\lIf{$|U|=1$}
{ $c=1$\;}}
$t:=2t$\;}
\Return{$U$}\;}
\caption{\label{al:agg} A bounding chain based on aggregation}
\label{algo_F}
\end{algorithm}

The lack of knowledge induced by the projection on $\cX$ forces us to consider \textit{all} the states $n\in\cS$ with the same projection $x$. This induces two main problems:
\begin{itemize}
 \item Even if the original system couples, we may never have $|U|=1$ in Algorithm \ref{algo_F}. 
 \item 
Even if at some time $-t$ we have only one value for the projected process (i.e. $|U|=1$), this is not necessarily the case for times $-s$, $-t\leq -s \leq 0$, as  
this single projected value at time $-t$ can correspond to many different states in $\cS$, as illustrated in Figure~\ref{fig:def_f}.
\end{itemize}

The first problem is similar as the case of EPSA algorithm when the envelope chain does not couple (i.e. the assumptions of Theorem \ref{th:env} do not hold). 
The second problem could not occur for EPSA algorithm (a singleton interval of an envelope chain contains only a single possible state of the initial chain, so only one possible trajectory from that time until time $0$). 
However, the two problems are of similar nature and the approaches described in Remark \ref{rem:epsa} remain valid. We will be interested in particular in performance bounds for increasing cost functions, 
described in Section \ref{ss:bounds}.

In addition to these problems, space $\cX$ can be too large to consider all the initial states $x\in\cX$. 
Our objective is twofold. First, we would like to be able to compute $\proj\left(\ffp(\Ux,a)\right)$, for a given state $x$ and event $a$, 
without considering all $n\in\Ux$ (for ATO systems with joint services, $|\Ux|$ grows exponentially with the number of item types). Second, we do not want to be forced to calculate, at each step, $\proj\left(\ffp(\Ux,a)\right)$ for \textit{all} the current states $x$ (even if we are able to calculate it easily for each state). To 
overcome this, we will combine the idea of aggregation with the method of envelopes developed in \cite{BGV08}. 

\subsection{Combining aggregation and envelopes} \label{subs:aggreg_env}

We will define a Markov chain $\Yinf$ (respectively $\Ysup$) that maps $x$ to the infimum (resp. supremum) of $f(x,a) = \proj\left(\ffp(\Ux,a)\right)$. 
Consider the following transition functions: for all $x\in\mathcal{X}$ and $a\in\cE$, 
\[
\underffp(x,a) \bydef \inf f(x,a), \quad \overffp(x,a) \bydef \sup f(x,a).
\]
Let $m,M\in\cX$ such that $m\peq M$. The envelope method in \cite{ashes} changes the subset $[m,M]$  into a new subset $[m',M']$  
(that depends on the event $a\in\cE$ and that usually involves only the transition function of one Markov chain). Here we consider two Markov chains $\Yinf$ and $\Ysup$ on the same space $\cX$, with the same set of events $\cE$, but with two different transition functions $\underffp$ and $\overffp$. Considering separately the envelopes of the infimum and the supremum chains does not necessarily sandwich the projected process (if either 
$\underffp$ or $\overffp$ is not monotone). 
We define the aggregated envelope transition function as follows:
For $m,M \in \cX$ such that $m \preceq M$, and $a\in\cE$, 
 \begin{eqnarray*}
\FFpp \left([m,M],a\right)  &\bydef&   \left[ \inf_{m\peq x \peq M} \underffp(x,a),\sup_{m\peq x \peq M} \overffp(x,a) \right] \\
&=& \left[\FFppinf\left([m,M],a\right),\FFppsup\left([m,M],a\right)\right].
 \end{eqnarray*}

In order to compare the projected process $Y=\proj(X)$ of the original chain $X$ to the lower and the upper envelopes of $\FFpp$, 
 we need the following notation. 
Assume the sequence of events $a_{-t+1},\dots,a_0\in\evtSet$ fixed.
For $(t,n)\in\mathbb{N}\times \cS$, $X{(-t,n)}$ stands for a realization of $X$ that starts from $n$ at time $-t$, while $X_{-s}{(-t,n)}$ denotes the value of this realization at time $-s$.
The next proposition shows that the chain with transition function $\FFpp$ is a bounding interval chain for the projected process $\proj (X)$.
\begin{proposition}
\label{lem:sandwiching_inf_sup}
Let $n\in\cS$ and $y, z \in\cX$ such that $y \peq \proj(n)\peq z$. Let $t\in\N$ and $a_{-t+1},\dots,a_0\in\evtSet$. Then we have for any $s\leq t$:
\begin{eqnarray*}
\label{eq:sup_inf}
\FFppinf \left([y,z],a_{-t+1\to -s}\right) &\peq& \proj(X_{-s}{(-t,n)}) \\
&\peq& \FFppsup \left([y,z],a_{-t+1\to -s}\right).
\end{eqnarray*}
\end{proposition}

\proof 
Let $x:=\proj(n)$. We prove the result by descending induction on $s$. First, for $s=t$ the result is trivial, as
$a_{-t+1\to -s} = \epsilon$ (empty word), so 
\begin{eqnarray*}
\FFppinf\left([y,z],\epsilon\right) = y \;\peq\; \proj(X_{-t}{(-t,n)})
\peq  z =  \FFppsup\left([y,z],\epsilon\right).
\end{eqnarray*}
Assume the result is true for some $s$, $s\leq t$. Let $n':=X_{-s}{(-t,n)}$, $x':=\proj(n')$, $y':=\FFppinf\left([y,z],a_{-t+1\to -s}\right)$ and $z':=\FFppsup\left([y,z],a_{-t+1\to -s}\right)$. Then we have that $y' \peq x' \peq z'$ by induction hypothesis. For $s-1$, we have:
\begin{eqnarray*}
\FFppinf\left([y',z'],a_{-s+1}\right) &\peq& \underffp(x',a_{-s+1}) \peq \proj(X_{-s+1}{(-t,n)}) \\
&\peq&  \overffp(x',a_{-s+1}) \peq  \FFppsup\left([y',z'],a_{-s+1}\right).
\end{eqnarray*}
By the definition of $y'$ and $z'$, $\FFppinf\left([y',z'],a_{-s+1}\right)=\FFppinf\left([y,z],a_{-t+1 \to -s+1}\right)$ and $\FFppsup\left([y',z'],a_{-s+1}\right)=\FFppsup\left([y,z],a_{-t+1 \to -s+1}\right)$, which gives the result for $s-1$ and ends the proof.
\endproof

The aggregated envelope method is summarized in Algorithm \ref{algo_A}.

\medskip

\begin{algorithm}[H]
\KwData{
\Iid events  $\left(a_{-t}\right)_{t\in \integer} \in \Evt^\N$.
}
\KwResult{Interval $[m^*,M^*]\subset\cX $ containing the projection of a state $n^*\in\cS$ distributed according to $\pi$.
}
  \Begin{
$t=1$; $c=0$\;
  \Repeat{$c=1$}
{ $m :=  \bot$ ($\in\cX$); $M :=  \top$ ($\in\cX$)\; 
 \For{$i=t-1$ {\em  \textbf{downto}} $0$}
{ $[m, M] :=\FFpp \left( [m, M],a_{-i}\right)$ \;
\lIf{$m=M$}
{ $c=1$\;}}
$t:=2t$\;}
$m^* := m$; $M^* := M$\;
\Return{$m^*$, $M^*$}\;}
\caption{Aggregated Envelope Perfect Sampling}
\label{algo_A}
\end{algorithm}

\subsection{Alternative stopping condition}
\label{ss:another}

Note that Algorithm \ref{algo_A} only gives an interval that contains the projection of a state distributed according to the stationary distribution. 
In Algorithm \ref{algoB}, we use a heuristic to relax the stopping condition:
Instead of stopping when the upper and lower envelopes meet ($m=M$), we stop when they meet at least once on each component between time $-t+1$ and time $0$. 
The intuition behind this is the following: For a chain which is componentwise coupling (see Definition \ref{def:cc} in Section \ref{subsubs:POS_slightly_general}), 
this condition is sufficient to insure that the chain has coupled. 
In the general case, this algorithm can be faster than Algorithm \ref{algo_A} (and still obtain bounds on the projection of a steady state). 
In addition, it is easier to bound the stopping time of Algorithm \ref{algoB}. However, as the chain has not necessarily coupled, there can be cases for which bounds obtained using Algorithm \ref{algoB} are loose. 

In the case of the ATO-POS system with joint returns, we are able to provide some bounds for the stopping time of Algorithm \ref{algoB} (see Theorem \ref{th:algo_B} in Section \ref{ss:POSgeneral}). Moreover, we will see that the results obtained for Algorithms \ref{algo_A} and \ref{algoB} are very close. 

\medskip

\begin{algorithm}[H]
\KwData{
\Iid events  $\left(a_{-t}\right)_{t\in \integer} \in \Evt^\N$.
}
\KwResult{Interval $[m^*,M^*]\subset\cX $ containing the projection of a state $n^*\in\cS$ distributed according to to the stationary distribution.} 
  \Begin{
$t:=1$; $c:= zeros(1,I)$\;
  \Repeat{$c=ones(1,I)$}
{ $m :=  \bot$ ($\in\cX$); $M :=  \top$ ($\in\cX$)\;
 \For{$i=t-1$ {\em  \textbf{downto}} $0$}
{ $[m, M] :=\FFpp \left( [m, M],a_{-i}\right)$ \;
\For{$j=1$ {\em  \textbf{to}} $I$}
{ \lIf{$m(j)=M(j)$}
{ $c(j):=1$\;}}}
$t:=2t$\;}
$m^* := m$; $M^* := M$\;
\Return{$m^*$, $M^*$}\;}
\caption{Modified stopping condition}
\label{algoB}
\end{algorithm}

\section{ATO systems with joint returns}
\label{sec:Aggreg_POS}

The case of ATO systems with joint returns is considerably more difficult, due to the need to track the information on which items were sent together to the customer. 
This leads to extremely large state space, not having a natural lattice structure (as will be shown in Section \ref{subs:POS_joint_services}), which makes it  
intractable for PSA or EPSA algorithms. We show in this section how to use the aggregated envelope method proposed in Section \ref{sec:method} to reduce the state space and overcome the lack of lattice structure. 

Even if the general idea is also valid for ATO-TOS, we focus only on ATO-POS for two reasons: unlike the TOS case, POS model does not have a product form solution, which makes it even more challenging. 
The second reason is that the monotonicity of arrivals in the POS case makes computations easier to some extent and simplifies the presentation  of the general ideas.  
The adaptation for the TOS case is discussed in Section \ref{sec:frc} (although TOS has product form, the estimation of the normalizing constant is still a hard problem).

We start by defining the model  and its event representation in Sections \ref{subs:POS_joint_services} and \ref{sec:POS_joint_services_inf}. 
In order to compute the aggregated envelope chain $\FFpp$ (Section \ref{app:inf}), we first analyze separately the supremum chain (Section \ref{app:sup}) and the infimum chain (Section \ref{subs:inf}).
In Section \ref{ss:POSgeneral}, we give a bound on the mean stopping time of Algorithm \ref{algoB}. 
When the service rate is high (Section \ref{sec:POS_joint_returns}), one can obtain perfect samples from the steady state distribution of the Markov chain, 
with an algorithm whose complexity is quadratic with respect to the total capacity $|C| = \sum_{i=1}^I C_i$.

\subsection{Model description} \label{subs:POS_joint_services}

We consider an Assemble-To-Order system with Partial Order Service (ATO-POS) and joint returns of items, also called 
the service tools model by \cite{VH09}. As in Section \ref{sec:POS}, demands for each subset $A\subset \I=\{1,...,I\}$ follow a Poisson process of rate $\lambda_A$. 
As before, we assume that the number $|\{A\subset \I:\lambda_A\neq 0\}|$ of subsets a customer can ask for is small 
(for instance linear with respect to the number $I$ of item types).
If some demanded items are not available, 
then the customer takes the available items (POS case). The available items return from the customer together (unlike the individual replenishment assumption in previous sections), 
and after an exponential time of rate $\mu$.

In terms of a network of $I$ queues, this means that we consider joint services: items that arrived in the system together (taken by one customer) will also leave the system together (returned from the customer). 
Therefore we need to keep track of which items arrived to the queues together. 
The system can be modeled as a continuous time Markov chain with 
state space:
\[\cS=\left\{(n_A)_{A\subseteq\I, A \neq \emptyset}\,:\,\forall A, n_A\geq 0 \;\&\; \forall i,\sum_{A:i\in A} n_A\leq C_i\right\} .\]
where $n_A$ is the number of subsets $A$ currently at the customers. 
Let $e_A$ denote the vector of $\cS$ whose coordinate $A$ is equal to $1$, and others are $0$ (this notation also stands for $\sum_{i\in A}e_i\in\cX$, but no confusion will be possible since the space will always be specified).
We define the projection $\proj$ on space $\cX=\{0,...,C_1\}\times \{0,...,C_2\} \times \cdots \times \{0,...,C_I\}$: 
\[\proj:\left\{\begin{array}{ccc} \cS & \longrightarrow & \mathcal{X} \\ n=(n_A)_{A\subset \I} & \longmapsto & x=\left(\sum_{A:i\in A} n_A\right)_{i\in \I} \\ \end{array} \right.\]
The total number of items of each type in queues (i.e. currently used by the customers) is given by a vector $x = (x_1, \ldots, x_I) \in\cX$, where $x_i$ is the number of items of type $i$. 
We will consider the product order $\leq$ on $\cX$.

Note that the only dimension of the state space $\cS$ is $2^I-1$. In addition, 
it is easy to see that $\cS$ equipped with the usual product order is not a lattice, as 
the supremum of two points can exceed the state space. For instance, let $I=2$ and set $c:=\min(C_1,C_2)$. Let $n,m\in\cS$ such that $n_{\{1,2\}}=c$, $n_{\{1\}}=n_{\{2\}}=0$ and $m_{\{1\}}=m_{\{2\}}=c$, $m_{\{1,2\}}=0$. 
Then $\proj(n)=\proj(m)=(c,c)$, yet the supremum $M=\sup(n,m)$ satisfies: $\proj(M)=(2c,2c)$. 
Therefore, PSA or EPSA cannot be applied directly using the product order on $\cS$: Even with monotone events, it is impossible to upper-bound any given subset of states by only a few extremal states.   

We have two different types of transitions. For each $n\in\cS$, and for each $A \subset \I$:
\begin{itemize}
 \item There is a demand for subset $A$, with rate $\lambda_A$. The new state is:
$
n +e_{A^{(n)}},
$
where 
\begin{equation}
\label{eq:An}
A^{(n)} = \left\{i\in A: (\proj(n))_i < C_i\right\}
\end{equation}
 denotes the items of set $A$ 
that are available in state $n$ and that are sent together to the customer. 
 \item If $n_A > 0$,  there is a joint service of $A$, with rate $\mu \cdot n_A$. The new state is $n-e_A$.
\end{itemize}

By a standard uniformization procedure, we can transform the above continuous time Markov chain to a discrete time Markov chain. 
The outgoing rate for each state is upper-bounded by
$
\Lambda:=\lambda + \mu |C|.
$
We take the uniformization constant to be equal to $\Lambda$. We now give a discrete event representation of the uniformized Markov chain.

\subsection{Discrete event representation} \label{sec:POS_joint_services_inf}

\paragraph{Arrivals.}
For any $A \subset \I$, $A \neq \emptyset$, let $d_{A}$ be the event of probability $\lambda_{A}/\Lambda$ that corresponds to a ``joint arrival to queues in $A$''.  
We give the transition function $\ffp$ of the Markov chain $X$ on $\cS$ for an arrival $d_A$, $ A\subset \I$, $A\neq \emptyset$, in state $n\in\cS$: 
\begin{eqnarray}
\label{eq:arrivals_X}
 \ffp(n,d_A)=n+e_{A^{(n)}},
\end{eqnarray}
where $A^{(n)}\subset A$ is defined in (\ref{eq:An}).

\paragraph{Services.}
Unfortunately, if service events are not well chosen, the supremum chain does not move with any service. Indeed, let us 
observe a fixed 
$i \in \I$. 
In order to have $(\overffp(x,a))_i = x_i - 1$, for some service event $a$, 
in \textit{all} states $n\in\cS$ such that $\proj(n)=x$, event $a$ must correspond to a service of some subset $A\subset\I$ that contains $i$ ($A$ can depend on $n$). Otherwise, if there is at least one state $n$ whose $i$-th queue is not served, then the $i$-th component of the supremum does not move.
This makes the definition of services a little tricky. In addition, the following representation is linear with respect to $|C|$.

Before defining service events, we need to define an ordering for the non-empty subsets of $\I$.
 For all $i\in\I$, we define $\left(A^i_k\right)_{0\leq k\leq 2^{I-i}-1}$ as an ordering of all the subsets of $\{i,...,I\}$ containing $i$ (subsets $\left(A^1_k\right)_k$ are those that contain 1, 
$\left(A^2_k\right)_k$ those that contain 2 but not 1, and so on). More precisely, let $i\in\I$ and $k\in\{0,1,...,2^{I-i}-1\}$. 
We set $k=k_1 \, ...\: k_{I-i}$ for the binary representation of $k$ ($k=\sum_{s=1}^{I-i} k_s 2^{I-i-s}$, where $k_s\in\{0,1\}$). Then the subset $A^i_k$ is by definition such that:
 \begin{itemize}
\item $i\in A_k^i$ ;
 \item for $s\in\{i+1,...,I\}$, $s\in A_k^i$ if and only if $k_{s-i}=0$.
\end{itemize}
For instance, if $I=4$, $A^1_0 = \{1,2,3,4\}$, $A^1_1 = \{1,2,3\}$, \ldots, $A^1_7 = \{1\}$, $A^2_0 = \{2,3,4\}$, \ldots, $A^4_0 = \{4\}$.
%
%

Now we can define services.
Let $r^i_j$, $1\leq i\leq I$, $1\leq j\leq C_i$, be independent events of probability $\mu/\Lambda$, such that for $n\in\cS$: 
\begin{itemize}
\item If $\sum_{\ell=0}^{2^{I-i}-1} n_{A^i_\ell} < j$, set $\ffp(n,r^i_j)=n$ ;
 \item If $\sum_{\ell=0}^{2^{I-i}-1} n_{A^i_\ell} \geq j$, let $k$ be the minimal element of $\{0,1,...,2^{I-i}-1\}$ such that ${\sum_{\ell=0}^{k-1} n_{A^i_\ell} < j} \leq \sum_{\ell=0}^{k} n_{A^i_\ell}$. 
Set $\ffp(n,r^i_j)=n-e_{A^i_k}$.
\end{itemize}
Note that the total number of events corresponding to services is exactly $|C|$. 
The next lemma shows that this definition of services agrees with the rates of the Markov chain defined in \ref{subs:POS_joint_services}.

\begin{lemma}
\label{lem:services}
 For any state $n\in\cS$ and any subset $A^i_k$, $1\leq i\leq I$, $0\leq k\leq 2^{I-i}-1$, with the above definition of services, the total probability to go from state $n$ to state $n-e_{A^i_k}$ is $\frac{\mu}{\Lambda}\cdot n_{A^i_k}$.
\end{lemma}

\proof 
Let $n \in \cS$ and $i \in \I$ be fixed. We assume that $n \not= 0$ (otherwise the result is obvious). 
First we explain the effect of services $(r^i_j)_j$ on $n$. For any $j\in\{1,...,C_i\}$, the service $r^i_j$ corresponds 
to a joint service of $A^i_k$, for some $k$ that depends on state $n$. More precisely, for $j=1$, service $r^i_1$ corresponds to a joint service of the subset $A^i_k$, where $k$ 
is the smallest integer such that $n_{A^i_k}>0$. We set $n'=n-e_{A^i_k}$. If $n'=0$, service $r^i_1$ is the only service that modifies $n$. 
Otherwise, for $j=2$, service $r^i_2$ corresponds to a joint service of the subset $A^i_{\ell}$, where $\ell$ is the smallest integer such that $n'_{A^i_\ell}=n_{A^i_\ell}-\mathds{1}_{\{\ell=k\}}>0$. 
Set $n''=n'-e_{A^i_{\ell}}$. If $n''= 0$, then $r^i_1$ and $r^i_2$ are the only services that change $n$, and if $n''\neq 0$, we continue as before.

The probability to go from state $n$ to $n-e_{A^i_k}$, for any $0 \leq k \leq 2^{I-i}-1$ is: 
$$
 \sum_{j=1}^{C_i}\Pb\left(\ffp(n,r^i_j)=n-e_{A^i_k}\right)  =  {\frac{\mu}{\Lambda}  \left|\left\{j\in\{1,\ldots,C_i\}:\sum_{\ell=0}^{k-1} n_{A^i_\ell} < j\leq \sum_{\ell=0}^{k} n_{A^i_\ell}\right\}\right|} 
{= \frac{n_{A^i_{k}}\mu }{\Lambda}.}
$$
\endproof

We want to compute the aggregated envelope transition function $\FFpp$ (defined in \ref{subs:aggreg_env}). 
For this, we first need to compute the supremum chain $\overffp$ and the infimum chain $\underffp$ (defined in \ref{subs:aggreg_env}). We begin by the supremum chain. 

\subsection{Supremum chain}
\label{app:sup}

The supremum chain was defined as follows in Section~\ref{subs:aggreg_env}. For all $a\in\cE$:
\begin{equation*}
\overffp(x,a)=\sup\left\{\proj\left(\ffp(n,a)\right):n\in\mathcal{S},\:\proj(n)=x\right\}.
\end{equation*} 

\begin{lemma}
\label{lem:comp_sup}
For all $x\in\cX$, set $\hat{x}_i\bydef\max\{x_i-(x_1+...+x_{i-1}),0\}$. 
Then the transition function of the supremum chain is, for $A\subset \I$, $A\neq\emptyset$, $i\in \I$, $1\leq j\leq C_i$, and $x\in\mathcal{X}$:
\begin{eqnarray*}
\left\{ \begin{array}{rcl} 
        \overffp(x,d_A) &=& x+\displaystyle\sum_{k\in A}\mathds{1}_{\{x_k<C_k\}}e_k = \fp(x,d_A), \\
	\overffp(x,r^i_j) &=&  x-\mathds{1}_{j\leq \hat{x}_i} e_{i}.
        \end{array}
\right.
\end{eqnarray*}
where $\fp$ is the transition function for the model ATO-POS with individual returns.
\end{lemma}

\proof 
{\em Arrivals.} Let $A\subset \I$, $A\neq\emptyset$, $x\in\cX$, and $n\in\cS$ such that $\proj(n)=x$. Applying $\proj$ to \eqref{eq:arrivals_X} gives that: 
$ \proj(\ffp(n,d_A))=x+\sum_{i\in A} \mathds{1}_{\{x_i <C_i \}} e_i$, 
hence the set $f(\{x\},d_A)=\left\{x+\sum_{i\in A} \mathds{1}_{\{x_i <C_i \}} e_i\right\}$ has only one element (see (\ref{eq:agf}) for the definition of $f$), and 
$\overffp(x,d_A)=x+\sum_{i\in A} \mathds{1}_{\{x_i <C_i \}} e_i=\underffp(x,d_A)$.

{\em Services.} Let $i\in \I$, $1\leq j\leq C_i$, and $x\in\mathcal{X}$. We first assume $j\leq \hat{x}_i$ (note that this implies $x_i\geq 1$). For $n$ such that $\proj(n)=x$, we will prove that
\begin{eqnarray}
\label{eq:comp_sup}
\left(\proj(\ffp(n,r^i_j))\right)_i=x_i-1,
\end{eqnarray}
which implies that $\overffp(x,r^i_j)\leq x-e_i$. In addition, if $n_{\{i\}}=x_i$, then $\proj\left(\ffp(n,r^i_j)\right)=x-e_i$. So $\overffp(x,r^i_j)= x-e_i$. 
Now we prove \eqref{eq:comp_sup}. Let $n\in\cS$ such that $\proj(n)=x$. Then 
\begin{eqnarray*}
 x_i=\sum_{A:i\in A} n_A=\sum_{i'=1}^i \sum_{k:i\in A^{i'}_k} n_{A^{i'}_k} \leq x_1+...+x_{i-1}+\sum_{k} n_{A^i_k}.
\end{eqnarray*}
where the second equality comes from the recursive definition of subsets $\left(A^i_k\right)_k$.
 Hence we have: $\sum_{k} n_{A^i_k} \geq \max\{x_i-(x_1+...+x_{i-1}),0\}=\hat{x}_i \geq j$, so we get the result by the definition of the event $r^i_j$.

We prove the second case now: assume $j > \hat{x}_i$. For $n$ such that $n_{\{k\}}=x_k$ for all $1\leq k\leq I$ (and $n_A=0$ if $|A|\geq 2$), we have: $\ffp(n,r^i_j)=n$, so $\overffp(x,r^i_j)= x$, which concludes the proof.

\subsubsection*{Monotonicity.} 
As a corollary of Lemma \ref{lem:comp_sup}, we have the following result (the proof for services follows from the fact that $x\leq y$ and $x_i=y_i$ imply $\hat{x}_i\geq\hat{y}_i$).

\begin{proposition}
\label{prop:mono_sup}
 Let $a\in\cE$ be any arrival or service. Under the usual product order on $\cX$, the event $a$ is monotone for the supremum chain in the ATO-POS system with joint services, \ie 
for all $x,y\in\cX$ such that $x\leq y$, we have that: $\overffp(x,a) \leq \overffp(y,a)$. 
\end{proposition}

\subsection{Infimum chain} \label{subs:inf}

The infimum chain was defined as follows in Section~\ref{subs:aggreg_env}, for all $x \in \cX, \; a\in\cE$:
\begin{equation*}
\underffp(x,a)=\inf\left\{\proj\left(\ffp(n,a)\right):n\in\cS,\:\proj(n)=x\right\}.
\end{equation*} 

For arrivals, we have for any $A\subset \I$, $A\neq\emptyset$, and $x\in\cX$ (see the proof of Lemma \ref{lem:comp_sup}): 
$$\underffp(x,d_A)=x+\sum_{i\in A} \mathds{1}_{\{x_i <C_i \}} e_i=\overffp(x,d_A).$$
The following lemma describes the infimum chain for services.
\begin{lemma}
\label{lem:comput_inf}
Let $p\in\I$. The $p$-th component of $\underffp(x,r^i_j)$ satisfies: 
\begin{itemize}
 \item If $p<i$, then $\left(\underffp(x,r^i_j)\right)_p=x_p$;
\item If $p=i$, then $\left(\underffp(x,r^i_j)\right)_p=x_i-\mathds{1}\{j\leq x_i\}$;
\item If $p>i$, then $\left(\underffp(x,r^i_j)\right)_p = x_p-\mathds{1}\left\{x_p > 0 \;\, \& \;\, j\leq \min\left(\sum_{i'=i+1}^p x_{i'},x_i\right)\right\}$.
\end{itemize}
\end{lemma}

\proof 

Since we consider the product order on $\cX$, the infimum can be computed componentwise, \ie for $p\in\I$:
$ \left(\underffp(x,r^i_j)\right)_p=\inf\left\{\left(\proj\left(\ffp(n,r^i_j)\right)\right)_p:n\in\cS,\:\proj(n)=x\right\}.
$

{\em Case $p<i$.} By the definition of $r^i_j$, there is no change on the $p$-th component, so $\left(\underffp(x,r^i_j)\right)_p=x_p$.

{\em Case $p=i$.} If $j> x_i$, then no $n$ such that $\left(\proj(n)\right)_i=x_i$ can be modified by $r^i_j$, so $\left(\underffp(x,r^i_j)\right)_i=x_i$. If $j\leq x_i$, let us choose $n$ such that $\proj(n)=x$ and $n_{\{i\}}=x_i$. Then $\left(\ffp(n,r^i_j)\right)_i=x_i-1$, so $\left(\underffp(x,r^i_j)\right)_i=x_i-1$. 

{\em Case $p>i$.} Clearly, a state $n\in\proj^{-1}(\{x\})$ can be modified by $r^i_j$ only if $x_p>0$. Moreover, 
its $p$-th component is modified by $r^i_j$ if and only if there exists $k\in\{0,1,...,2^{I-i}-1\}$ such that 
\begin{eqnarray}
\label{eq:move}
 p\in A^i_k \quad\textrm{and}\quad \sum_{\ell=0}^{k-1} n_{A^i_\ell} < j \leq \sum_{\ell=0}^{k} n_{A^i_\ell}.
\end{eqnarray}
In particular, the condition $j \leq \sum_{\ell=0}^{2^{I-i}-1} n_{A^i_\ell}$ (and so $j\leq x_i$) is necessary for $n$ to change. 

Assume $j > \sum_{i'=i+1}^p x_{i'}$. We will show that no $n$ such that $\proj(n)=x$ can change on the $p$-th component by $r^i_j$. Indeed, 
 due to the recursive definition of subsets $(A^i_k)_k$,
we have that $A^i_k$ cannot contain $p$ if $k\geq 2^{I-i-1}-2^{I-i+1-p}$, and, for $k=2^{I-i-1}-2^{I-i+1-p}-1$, we have:
\begin{eqnarray*}
 \sum_{\ell=0}^{k} n_{A^i_\ell} =\sum_{i'=i+1}^p \sum_{\begin{array}{c} \scriptstyle \ell:i'\in A^i_\ell \\ \scriptstyle 1<i''<i' \Rightarrow i'' \notin A^i_\ell  \\ \end{array} } n_{A^i_\ell}
 \quad \leq \sum_{i'=i+1}^p x_{i'}
< j.
\end{eqnarray*}
So \eqref{eq:move} cannot hold, and $n$ cannot change on the $p$-th component.

Assume now $x_p>0$ and  $j \leq \min\left(\sum_{i'=i+1}^p x_{i'},x_i\right)$. We give the construction of $n\in\cS$ such that $\proj(n)=x$ and $\left(\ffp(n,r^i_j)\right)_p=x_p-1$: 

\medskip

\begin{center}
{
\begin{tabular}{|l|}
 \hline
\textbf{Construction of $n$} \\
\hline
$m:=\min(x_i,x_p)$ ;  \\
$V:=0$ ; $a:=0$ ; \\
\textbf{for} $i'=i+1$ \textbf{to} $p-1$ \textbf{do} \\
\qquad $a:=\min(x_{i'},j-1-V,x_i-m-V)$ ; \\
\qquad $n_{\{i,i'\}}:=a$ ; \: $n_{\{i'\}}:=x_{i'}-a$ ; \\
\qquad $V:=V+a$ ; \\
\textbf{end} \\
$n_{\{i,p\}}:=m$ ; \: $n_{\{p\}}:=x_p-m$ ; \: $n_{\{i\}}:=x_i-V-m$ ; \\
$n_{\{i'\}}=x_{i'}$ for $i'<i$ or $i'>p$, and $n_A=0$ for other subsets $A$ ; \\
\hline
\end{tabular}
}
\end{center}

\medskip 

Clearly, $\proj(n)=x$. Note also that our hypotheses imply that $m\geq 1$.
Let $k\in\{0,1,...,2^{I-i}-1\}$ such that $A^i_k=\{i,p\}$. We will show that:
$
 \sum_{\ell=0}^{k-1} n_{A^i_\ell} < j \leq \sum_{\ell=0}^{k} n_{A^i_\ell},
$
which will end the whole proof (cf. \eqref{eq:move}). Let $\tilde{V}$ be the value of $V$ at the end of the loop. Then we have that $\sum_{\ell=0}^{k-1} n_{A^i_\ell}=\tilde{V}$ and $\sum_{\ell=0}^{k} n_{A^i_\ell}=\tilde{V}+m$, 
so we are left to prove that:
\begin{eqnarray}
\label{eq:S}
 \tilde{V} < j \leq \tilde{V}+m
\end{eqnarray}
We study more precisely what happens during the loop. If there exists a step during the loop such that $a=j-1-V$ (resp. $a=x_i-m-V$), 
then $a=0$ in all the next steps and, at the end of the loop, $\tilde{V}=j-1$ (resp. $\tilde{V}=x_i-m$). In these two cases, \eqref{eq:S} follows easily. 
The remaining case is the one such that, at each step of the loop, $a=x_{i'} \leq \min(j-1-V,x_i-m-V)$. 
In particular, considering the last step of the loop, the previous inequality gives that $\tilde{V} \leq j-1$. Moreover, $\tilde{V}=\sum_{i'=i+1}^{p-1} x_{i'}$ in this case, so the fact that $j \leq \min\left(\sum_{i'=i+1}^p x_{i'},x_i\right)$ gives that $j\leq \tilde{V}+m$ (considering the two possible cases for the value of $m$).
\endproof

\subsubsection*{Non-monotonicity.}

Unfortunately, services are not monotone for the infimum chain. Indeed, let $I=2$, $x=(0,1)$, and $y=(1,1)$. Then $x\leq y$, yet 
$
 \underffp(x,r^1_1) = (0,1) \geq (0,0) = \underffp(y,r^1_1).
$

\subsection{The aggregated envelope chain $\FFpp$}
\label{app:inf}

We use Sections \ref{app:sup} and \ref{subs:inf} to compute the aggregated envelope chain $\FFpp$ (defined in Section~\ref{subs:aggreg_env}).
Let $m,M\in\cX$, $m\leq M$. Our goal is to compute, for all $a\in\cE$:
 \begin{eqnarray*}
\FFpp \left([m,M],a\right)  &=&   \left[ \inf_{m\leq x \leq M} \underffp(x,a),\sup_{m\leq x \leq M} \overffp(x,a) \right] \\
&=& \left[\FFppinf\left([m,M],a\right),\FFppsup\left([m,M],a\right)\right].
 \end{eqnarray*}

\textit{Arrivals.}  Let $A\subset \I$, $A\neq \emptyset$. Then $\FFpp \left([m,M],d_A\right)=\left[\fp(m,d_A),\fp(M,d_A)\right]$ since 
$ \underffp(x,d_A)=\overffp(x,d_A)=x+\sum_{k\in A}\mathds{1}_{\{x_k<C_k\}}e_k = \fp(x,d_A) $ for all $x\in\cX$. 
In addition, $\fp$ is monotone due to Proposition \ref{prop:ATOpos_monotone}, which gives the result for $\FFpp$.

\textit{Services.} Let $i\in\I$ and $j\in\{1,2,...,C_i\}$. 
By Proposition \ref{prop:mono_sup}, $\overffp$ is monotone, thus 
$\FFppsup \left([m,M],r^i_j\right)=\overffp(M,r^i_j) = M-\mathds{1}_{j\leq \hat{M}_i} e_{i}.$
In order to compute $\FFppinf \left([m,M],r^i_j\right)$, we need to compute the lower envelope of the infimum chain (Lemma~\ref{lem:comput_env_inf}).

\begin{lemma}
\label{lem:comput_env_inf}
Let $m$, $M\in\mathcal{X}$ such that $m\leq M$. Let $i\in\I$ and $j\in\{1,2,...,C_i\}$. Set
$
m' = \FFppinf\left([m,M],r^i_j\right).
$
We compute each component $p$ of $m'$, for $p\in\I$, and we can distinguish three cases: 
\begin{itemize}
 \item If $p<i$, then $m'_p=m_p$; 
\item If $p=i$, then $m'_p=m_i-\mathds{1}\{j\leq m_i\}$;
\item If $p>i$, then $m'_p=m_p-\mathds{1}\left\{m_p > 0 \;\, \& \;\, j\leq \min\left(\sum_{i'=i+1}^{p-1} M_{i'} + m_p , M_i\right)\right\}$. 
\end{itemize}
\end{lemma}

\proof 
{\em Case $p<i$.} By the definition of $r^i_j$, there is no change on the $p$-th component, so $m'_p=m_p$.

{\em Case $p=i$.} Let $x,y\in\mathcal{X}$ such that $x\leq y$. Then:
$
 \left(\underffp(x,r^i_j)\right)_p=x_i-\mathds{1}\{j\leq x_i\} \leq y_i-\mathds{1}\{j\leq y_i\}=\left(\underffp(y,r^i_j)\right)_p.
$
Indeed, if $x_i<y_i$, this is obvious since $y_i$ looses at most 1. Otherwise, $x_i=y_i$ and $x$ and $y$ move together. 

{\em Case $p>i$.} If $j>M_i$, then for all $x\in\mathcal{X}$ such that $m\leq x\leq M$, we have that $j>x_i$, so $x$ cannot move and $m'_p=m_p$. Assume $j\leq M_i$ from now on. 
If $m_p=0$, then $m'_p=0$, so we assume $m_p>0$. We distinguish two cases:
\begin{itemize}
\item
Assume $j > \sum_{i'=i+1}^{p-1} M_{i'} + m_p$. Let $x\in\mathcal{X}$ such that $m\leq x\leq M$. If $x_p>m_p$, then $\left(\underffp(x,r^i_j)\right)_p\geq m_p$. Otherwise, $x_p=m_p$, and so $j > \sum_{i'=i+1}^{p-1} M_{i'} + m_p \geq \sum_{i'=i+1}^{p} x_{i'} $. Hence $\left(\underffp(x,r^i_j)\right)_p=x_p$. Finally this leads to $m'_p=m_p$.
\item
Assume $j \leq \sum_{i'=i+1}^{p-1} M_{i'} + m_p$. We have that $m'_p \geq m_p-1$, and we will define $x\in\mathcal{X}$ such that $m\leq x\leq M$ and $\left(\underffp(x,r^i_j)\right)_p=m_p-1$. This will then prove that $m'_p=m_p-1$. Set $x_p:=m_p$ and $x_{i'}:=M_{i'}$ for all $i'\in\{1,2,...,I\} \setminus \{p\}$. Thus $j\leq \min\left(\sum_{i'=i+1}^p x_{i'},x_i\right)$ and  $\left(\underffp(x,r^i_j)\right)_p=m_p-1$.
\end{itemize}
\endproof

Hence we gave the computation of the aggregated envelope chain $\FFpp$, which is necessary to use Algorithms \ref{algo_A} (in Section \ref{subs:aggreg_env}) and \ref{algoB} (in Section \ref{ss:another}). 

In order to give a bound on the complexity of Algorithm \ref{algo_A} or \ref{algoB}, we have to take into account:
\begin{itemize}
 \item the random choice of events,
\item the computation of the aggregated envelope transition function,
\item the stopping time of Algorithm \ref{algo_A} or \ref{algoB}.
\end{itemize}
 The number of services is linear with respect to $|C|$. In addition, the number of arrivals is equal to the number of subsets $A$ such that $\lambda_A\neq 0$: 
 If this number is of order $O(I)$, using alias method by \cite{Walker77}, the choice of events can be done in a constant time. 
 In addition, Lemma \ref{lem:comput_env_inf} shows that the computation of the aggregated envelope transition function can be done in a linear time with respect to $I$.
 We next study the stopping time of Algorithm \ref{algoB}.

\subsection{Bound for the stopping time of Algorithm \ref{algoB}} \label{ss:POSgeneral}

We suppose there exist two subsets $\I_0$ and $\I_C\subset\I$, $\I=\I_0 \cup \I_C$, such that:
\begin{itemize}
\item[(i)] $  \mu > \sum_{i\in \I_0} \lambda_i $,
\item[(ii)] $ \delta_p \bydef \lambda_p - \mu\left(\sum_{i=1}^p C_{i}-1\right)>0 $ for all $p\in \I_C$.
\end{itemize}
Without loss of generality, we can change the numbering of queues such that:
\begin{itemize}
 \item[(iii)] $\Big( \: i\in \I_0 \: \textrm{ and } \:j\in \I_C \: \Big) \; \Longrightarrow \; i\leq j$.
\end{itemize}

In the case where assumptions \textit{(i)} and \textit{(ii)} are not satisfied, we are not able to give a bound on the stopping time of Algorithm \ref{algoB} (using the same method for the proof). 
The last condition comes from the expression of $\FFppsup \left([m,M],r^i_k\right)= M-\mathds{1}_{k\leq \hat{M}_i} e_{i}$, where $\hat{M}_i=\max\{M_i-(M_1+...+M_{i-1}),0\}$. 
Indeed we want the upper process to reach $0$ on components $i\in\I_0$, and the lower process to reach $C_j$ on components $j\in\I_C$. 
Yet the upper process cannot decrease on the $i$-th coordinate (due to service $r^i_k$, for any $k$) if there exists $j<i$ such that $M_j$ is greater than $M_i$ 
(in that case, $\hat{M}_i=0$ and $\FFppsup \left([m,M],r^i_k\right)= M$).
\begin{theorem}
\label{th:algo_B}
 Assume conditions $(i)$ to $(iii)$ hold. Then we can bound the time $\tau_{\text{Alg}}$ for which all components couple at least once by:
$$
 \Eb\left[\tau_{\text{Alg}}\right]\leq \frac{\Lambda}{\mu-\sum_{i\in \I_0} \lambda_i} \sum_{i\in \I_0} C_i + \sum_{p \in \I_C} \frac{\Lambda}{\delta_p} C_p.
$$
\end{theorem}

The proof is given in Appendix \ref{app:thfive}, and requires Appendix \ref{app:Foster}, in which we give a bound on the mean hitting to zero for the supremum chain $\Ysup$.

We discuss the complexity of Algorithm \ref{algoB} (with respect to $|C|$) in the case where $\I_0=\I$ (high service rate case): Since $\Lambda=\lambda + \mu|C|$, 
we have that the mean stopping time of Algorithm  \ref{algoB} is of order $O\left(|C|^2\right)$ (using Theorem \ref{th:algo_B}). 
Due to the discussion at the end of Subsection \ref{app:inf}, the complexity of Algorithm \ref{algoB} is thus quadratic with respect to $|C|$ (for high service rate).

Recall this is only a stopping criterion: The chain does not necessarily couple if we use Algorithm \ref{algoB}, and it is interesting to compare it with Algorithm \ref{algo_A}, 
for which the stopping time is greater than or equal to the coupling time of the chain.

\begin{figure}
\begin{center}
\includegraphics[width=0.51\textwidth]{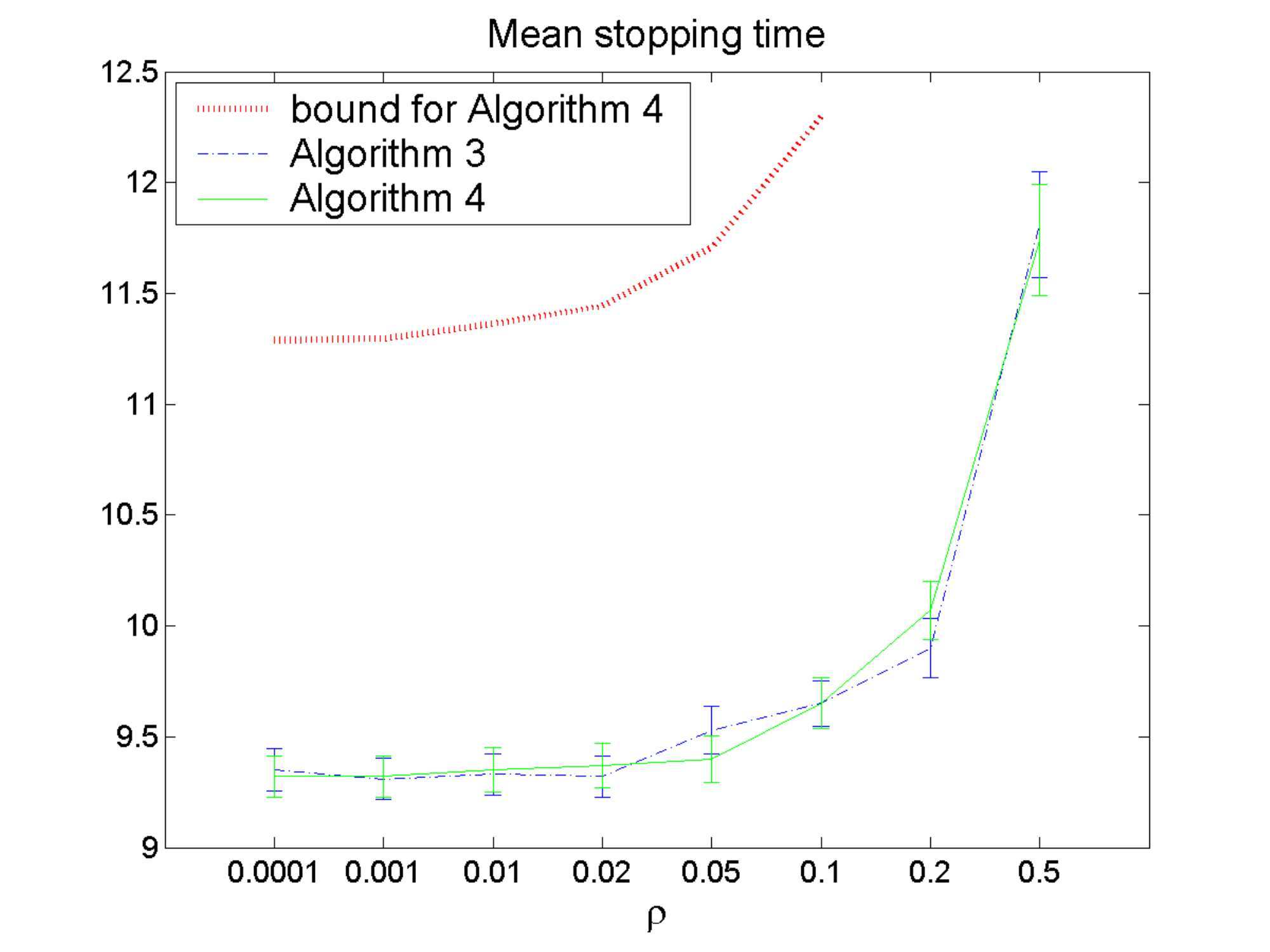}
\hspace{-2em}
\includegraphics[width=0.51\textwidth]{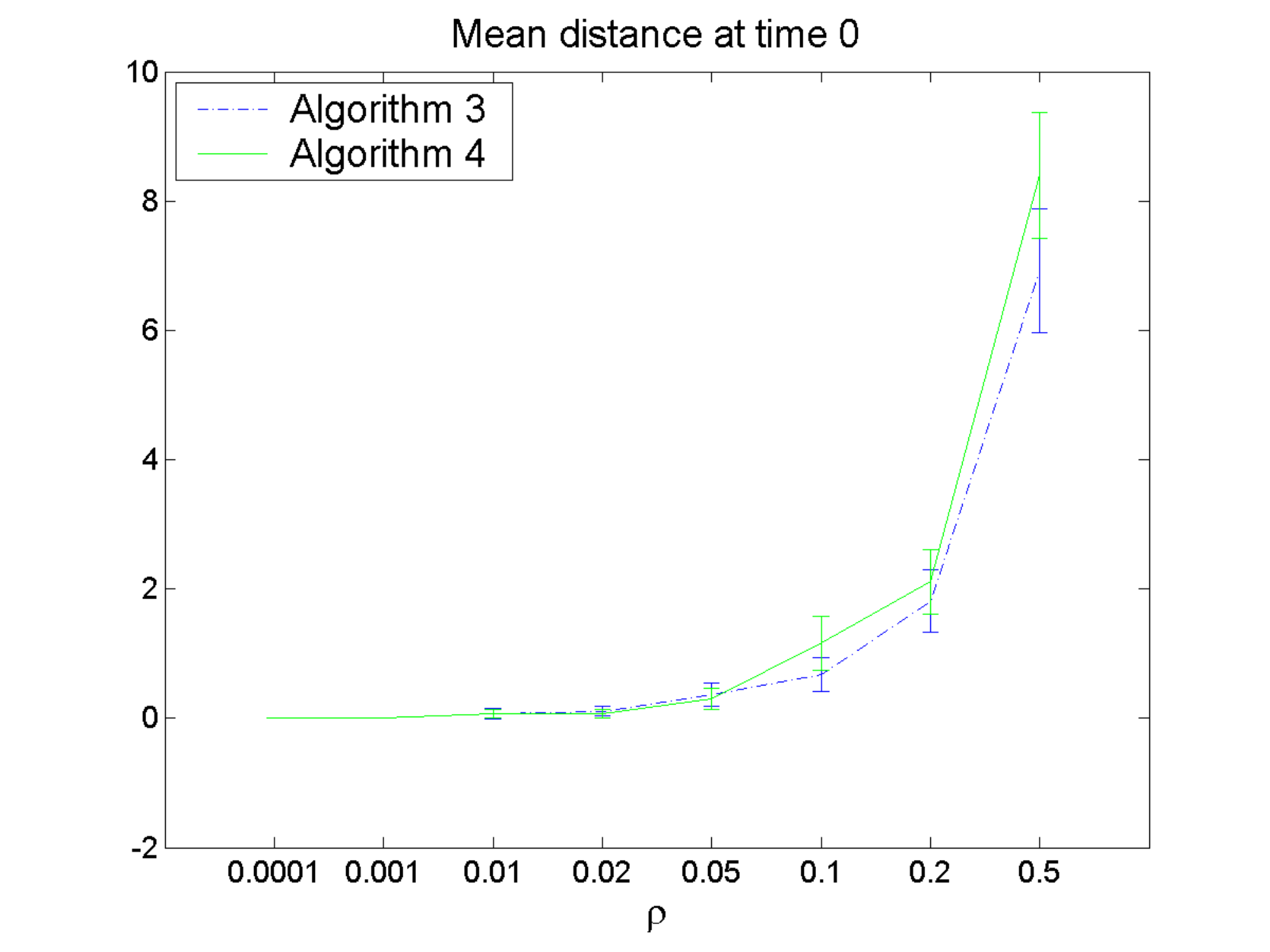}  
\caption{\label{fig:RJ} On the left: Stopping times for Algorithms \ref{algo_A} and \ref{algoB} (ATO-POS with joint returns), and the upper bound for Algorithm \ref{algoB} (we display $\log_2(T)$ where $T$ is the mean stopping time), together with the $95\%$ confidence intervals.
On the right: Mean distance (in $1$-norm) between upper and lower bounding states at time $0$ for Algorithms \ref{algo_A} and \ref{algoB}, together with the $95\%$ confidence intervals.}
\end{center}
\end{figure}

In Figure \ref{fig:RJ} on the left, we give stopping times for Algorithms \ref{algo_A} and \ref{algoB} (ATO-POS with joint returns), for the following parameters: $I=5$, $C_i=10, \forall i$, $\lambda_A = \frac{1}{2^{|A|-1}}$, and $\mu_{i}(x_i) = \mu x_i$, with $\rho = \frac{\lambda_i}{\mu}$. The size of the sample is $N=100$. We can observe that the mean stopping times of both algorithms are very close. In Figure \ref{fig:RJ} on the right, we provide mean distance between upper and lower bounding states at time $0$ using $1$-norm, \ie $\sum_{i \in \I} (M^*_i - m^*_i)$.

\subsection{Sampling the original (non aggregated) chain}
\label{sec:POS_joint_returns}

Instead of having an interval that contains the projection of a state distributed according to the stationary distribution, we can even sample exactly the steady state distribution. 
However, the algorithm we provide here is efficient only in the high service rate case. We prove in Theorem \ref{th:HSRCase} that its stopping time is quadratic with respect to the total capacity $|C|$.

The supremum chain $\Ysup$ is monotone (Proposition \ref{prop:mono_sup}, in Section \ref{app:sup}). This gives directly the computation for the upper envelope of $\FFpp$: 
for any $m,M\in\cX,a\in\cE$, $\FFppsup\left([m,M],a\right)=\overffp(M,a)$. This result and Lemma \ref{lem:sandwiching_inf_sup} (in Section \ref{subs:aggreg_env}) 
give that the projected chain $\proj(X)$ is between zero and the supremum chain $\Ysup$. If service rate is high, we can 
wait until the supremum chain $\Ysup$ hits zero. 
The main advantage is that this provides some solution for the possible decoupling of the system: When $\Ysup$ reaches zero, the projected process is also in state zero, 
and thus the \textit{only} possible state for the original Markov chain $X$ is also zero. Thus the original chain $X$ couples in $\cS$. 
Hence, we can use coupling from the past for $\Ysup$ until we find a time $-t$ such that $\Ysup_{-t}=0$, and then, from time $-t$ to time $0$, simulate the only trajectory of $X$ starting from state zero (with the same events). 
This algorithm has quadratic complexity in $|C|$, when $\mu>\sum_i \lambda_i$ (due to Lemma \ref{lem:Foster_sup} in Appendix \ref{app:Foster}).

We have the following bound for the mean coupling time of the original chain $X$, which is quadratic with respect to $|C|$:
\begin{theorem}
\label{th:HSRCase}
 Let $\tau$ be the coupling time of the original chain $X$ defined on $\cS$, and corresponding to the ATO-POS system with joint services. Assume $\mu>\sum_i \lambda_i$. Then we have:
\begin{equation*}
 \Eb[\tau] \leq \frac{\Lambda}{\mu-\sum_i \lambda_i} |C|.
\end{equation*}
\end{theorem}

\proof 
This is a consequence of Lemma \ref{lem:Foster_sup} in Appendix \ref{app:Foster} 
and the fact that the \textit{only} possible state for the chain $X$ when $\proj(X)=0$ is the state $0$: hence if the projected process reaches zero, the original system defined on $\cS$ also reaches zero, and so couples.
\endproof


\section{Further remarks and conclusions}
\label{sec:frc}

The contribution of our paper is twofold: 
\begin{itemize}
 \item We give perfect sampling algorithms for ATO models with individual and joint replenishments.
 \item We introduce a new method for perfect sampling, based on aggregation and bounding chains. We apply our new method to the ATO-POS case with joint  replenishments, for which the known perfect sampling techniques 
cannot be applied.
\end{itemize}

ATO systems with joint services have extremely large state space - its only dimension
is $2^{I}-1$ - and up to our best knowledge, there is no known efficient solution technique in the literature (in particular for the POS case that does not have a product form solution). 
Thus our new perfect sampling method can be of great interest to evaluate their performance, 
as well as 
in the optimization algorithms for capacity dimensioning. 
In most applications, the lost probability is demanded to be very low, thus the conditions in Section \ref{sec:POS_joint_returns} seem to be natural, under which we can obtain the exact samples of the stationary distribution. 

The results in Section~\ref{sec:Aggreg_POS} can be extended to the TOS case. The arrivals for the TOS case are no longer monotone, so 
the approach used in Section \ref{subs:inf} to compute the aggregated envelope chain for the infimum chain has to be used also for the supremum chain.
Note that the services are the same as for the POS model. 


Up to our knowledge, this is the first time that the aggregation of Markov chains is combined with perfect sampling technique to avoid state space explosion problems. 
This direction sounds promising for various applications.


\appendix

\section{Hitting time to zero for the ATO-POS system with individual state-dependent services (Section  \ref{subsubs:POS_slightly_general}).}
\label{sec:AppB}

We consider a slightly more general model than the one in Section \ref{sec:POS}: here we allow the service rate $\mu_i$ to depend on the whole state $x$ (not only on $x_i$), provided that services are monotone. This more general setting will be useful to prove Lemma \ref{lem:Foster_sup} (in Appendix \ref{app:Foster}).

More precisely, 
for each $1\leq i\leq I$, let $\mu_i:\cX \rightarrow [0,1]$, and set $\nu_i:=\max_{x\in\cX} \mu_i(x)$. Set $\Lambda:=\lambda+\sum_i \nu_i$ for the uniformization constant.
Arrival events are as before. The set $\varSigma$ of service events satisfies, for all $x,y\in\cX$ and $s\in\varSigma$:
\begin{eqnarray*}
 \: x\leq y \; \Rightarrow \; \fp(x,s)\leq \fp(y,s)\:  .
\end{eqnarray*}
 
Let $\tau_0^P$ be the mean hitting time that this ATO-POS system reaches the state $0$ in $\cX$. Previous assumption is satisfied by the model defined in \ref{subs:POS_state_dependence}, so the following proposition also holds for the initial model. 
\begin{proposition}
\label{prop:POS_zero}
 Let $\delta:=\min_{x\neq 0} \sum_i \mu_i(x)-\sum_i \lambda_i$ and assume $\delta>0$. Then we have:
$
 \Eb[\tau_0^P] \leq \frac{\Lambda}{\delta} \sum_{i=1}^I C_i.
$
\end{proposition}

\proof 
  Let $P$ be the transition matrix of the ATO-POS model described above. To prove this bound, we will use Lemma \ref{Foster}, with $S=\cX$, $U=\{0\}$ and, for all $z\in\cX$, $h(z):=\sum_{i=1}^I z_i$. Then we have, for all $y\in\cX\setminus \{0\}$:
\begin{eqnarray*}
 \sum_{z\in\cX} P(y,z) h(z) &=& \sum_A \frac{\lambda_A}{\Lambda} h\left(y+\sum_{i\in A}\mathds{1}\{y_i < C_i \}\: e_i\right)  +\sum_i \frac{\mu_i(y)}{\Lambda} h(y-e_i) + \sum_i \frac{\nu_i-\mu_i(y)}{\Lambda}h(y) \\
&=& \sum_A \frac{\lambda_A}{\Lambda} \left(\sum_i y_i + \sum_{i\in A} \mathds{1}\{y_i < C_i \} \right) +\sum_i \frac{\mu_i(y)}{\Lambda} \left( \sum_j y_j -1\right) + \sum_i \frac{\nu_i-\mu_i(y)}{\Lambda}\sum_j y_j \\
&=& \left( \sum_A \frac{\lambda_A}{\Lambda} +\sum_i \frac{\nu_i}{\Lambda} \right) h(y) + \sum_A \sum_{i\in A} \mathds{1}\{y_i < C_i \}  \frac{\lambda_A}{\Lambda} -\sum_i \frac{\mu_i(y)}{\Lambda} \\
&\leq & h(y) + \sum_A |A| \frac{\lambda_A}{\Lambda} -\sum_i \frac{\mu_i(y)}{\Lambda} = h(y) + \sum_i \frac{\lambda_i}{\Lambda} -\sum_i \frac{\mu_i(y)}{\Lambda} \leq h(y) -\frac{\delta}{\Lambda}.
\end{eqnarray*}
Hence the condition \eqref{eq:Foster_cond} of Lemma \ref{Foster} is proved, and we can apply \eqref{eq:Foster_ccl} with $y=(C_1,...,C_I)$. It follows that
$
 \Eb[\tau_0]=\Eb_y[\tau_{\{0\}}] \leq \frac{\Lambda h(y)}{\delta}=\frac{\Lambda}{\delta} \sum_{i=1}^I C_i
$
since the time for the system to hit zero is equal to the time for state $y=(C_1,...,C_I)$ to hit zero, due to monotonicity.
\endproof

\section{Proof of Proposition \ref{prop:I3} (Non-monotonicity of the ATO-TOS system with individual replenishments, Section  \ref{subs:envelopeTOS})}
\label{sec:AppI3}

Let $C=(C_1,...,C_I)$. Let $\preceq$ be a partial order such that $\ft$ is monotone for $\preceq$ (\ie all the events of $\cE$ are monotone). We proceed in three steps:
\begin{description}
 \item \textit{Step 1.} Let $x=(x_1,...,x_I)\in\cX$, $x\neq C$. Then we cannot have $x\preceq C$.
\item \textit{Step 2.} $C$ is not comparable to any other state.
\item \textit{Step 3.} Let $x\in\cX$. Then $x$ is not comparable to any other state.
\end{description}

\textit{Proof of step 1.} Assume by contradiction that $x\preceq C$. Since $x\neq C$, there exists $1\leq i\leq I$ such that $x_i<C_i$. Without loss of generality, we can assume that $i=1$. 
We first show that $v_1=(C_1-1,C_2,...,C_I)$ also verifies $v_1\preceq C$. Since $x\preceq C$ and $\ft$ is monotone for $\preceq$, 
it is enough to find a finite sequence of events that moves $x$ to $v_1$, and that let $C$ unchanged: For instance, we apply $C_1-1-x_1$ times event $d_{\{1\}}$  and then, 
for each $2\leq j\leq I$, we apply $C_j-x_j$ times event $d_{\{j\}}$. Hence $v_1=(C_1-1,C_2,...,C_I)\preceq C$.
Yet the combination $s^{(1)}_2$ followed by $d_{\{1,2\}}$ moves $v_1$ to $C$, and $C$ to $v_2=(C_1,C_2-1,...,C_I)$. Hence, using again the monotonicity of $\ft$, this leads to: $C\preceq {v_2=(C_1,C_2-1,...,C_I)}$.
With the same trick ($s^{(1)}_3$ followed by $d_{\{1,3\}}$), we have that $C\preceq {v_3=(C_1,C_2,C_3-1,C_4,...,C_I)}$.
 Yet, starting from $C\preceq v_2=(C_1,C_2-1,...,C_I)$ and using $s^{(1)}_3$ followed by $d_{\{2,3\}}$ leads to $v_3\preceq C$. Hence, $v_3= C$, which is a contradiction.

\textit{Proof of step 2.} Let $x\neq C$. The same argument as for Step 1 also works, starting from $C\preceq x$ instead of $x\preceq C$. Hence we cannot have $x\preceq C$ nor $C\preceq x$, 
which proves that $C$ is not comparable to any other state.

\textit{Proof of step 3.} Let $y=(y_1,...,y_I)\neq x$: there exists $1\leq i\leq I$ such that $x_i\neq y_i$. We assume by contradiction that either $y\preceq x$ or $x\preceq y$. 
Without loss of generality, we can assume that $y\preceq x$. 
Then $x_i>y_i$. We use the following sequence of events: For each $1\leq j\leq I$, we apply $C_j-x_j$ times the event $d_{\{j\}}$. 
This sequence moves $x$ to $C$ and $y$ to $y'$, where $y'$ is defined by: $y'_j=C_j-(x_j-y_j)$ for all $1\leq j\leq I$. In particular, $y'_i<C_i$, so $y'\neq C$. 
In addition, the fact that $y\preceq x$ and the monotonicity of $\ft$ for $\preceq$ imply that $y'\preceq C$, which is a contradiction.

\section{Hitting time to zero for the supremum chain (ATO-POS system with joint returns, Section  \ref{app:sup})}
\label{app:Foster}

The next lemma gives the mean hitting time to zero for the supremum chain, using the results of Appendix \ref{sec:AppB}  and Section \ref{app:sup}. 
It is used in the proofs of Theorem \ref{th:algo_B} (Section \ref{ss:POSgeneral}) and Theorem \ref{th:HSRCase} (Section \ref{sec:POS_joint_returns}).

\begin{lemma}
 \label{lem:Foster_sup}
Let $\bar{\tau}$ be the time that the supremum chain, starting from $C=(C_1,...,C_I)$, reaches the state $0$ in $\cX$. 
Assume $\mu>\sum_i \lambda_i$. Then we have:
$
 \Eb[\bar{\tau}] \leq \frac{\Lambda}{\mu-\sum_i \lambda_i} |C|.
$
\end{lemma}

\proof 
We apply Proposition \ref{prop:POS_zero} (Appendix \ref{sec:AppB}). Lemma \ref{lem:comp_sup}  and Proposition \ref{prop:mono_sup} show that the supremum chain is a particular case of the model presented in Appendix \ref{sec:AppB} 
(with $\mu_i(x)=\mu  \hat{x}_i$, $1\leq i\leq I$). The fact that $\delta=\min_{x\neq 0} \sum_i \mu_i(x)-\sum_i \lambda_i$ is positive comes from the fact that 
$\min_{x\neq 0} \sum_i \mu_i(x)=\mu\cdot \min_{x\neq 0} \sum_i \hat{x}_i=\mu$ and the hypothesis $\mu>\sum_i \lambda_i$.
\endproof

\section{Proof of Theorem \ref{th:algo_B} (Stopping time of Algorithm \ref{algoB}, Section \ref{ss:POSgeneral})}
\label{app:thfive}

The proof of Theorem \ref{th:algo_B} is based on Lemma \ref{lem:Foster_sup} (Appendix \ref{app:Foster}) and the following lemma, that gives a bound on the mean hitting time of $C_p$ for the $p$-th component of the infimum: 
\begin{lemma}
\label{lem:inf_capacity}
 Let $p\in\I$, and assume $\delta_p:=\lambda_p-\mu\left(\sum_{i=1}^p C_i -1\right)$ is positive. 
 Let $\tau^{(p)}$ be the time for the $p$-th component of $\FFppinf$ to hit $C_p$ (starting from $0$).
Then:
$
 \Eb\left[\tau^{(p)}\right]\leq \frac{\Lambda}{\delta_p} C_p.
$
\end{lemma}

\proof[Proof of Lemma \ref{lem:inf_capacity}.]
 In order to give a bound on $\tau^{(p)}$, we have to face two problems: $\FFppinf$ is not the transition function of a Markov chain (it depends on $\FFppsup$), and its projection on the $p$-th component depends on its whole state (see Lemma \ref{lem:comput_env_inf}). That is why we introduce a new Markov chain, defined on $\cX$, whose transition function $\Inf$ is a lower bound for $\FFppinf$, and whose projection on the $p$-th component is also a Markov chain. Indeed, for $m\in\cX$, $i\in\I$, $j\in\{1,...,C_i\}$, we define, for all $p\in\I$ (setting $m':=\Inf(m,r^i_j)$):
\begin{itemize}
 \item If $p<i$, then $m'_p=m_p$;
\item If $p=i$, then $m'_p=m_i-\mathds{1}\{j\leq m_i\}$; 
\item If $p>i$, then $m'_p=m_p-\mathds{1}\left\{m_p > 0 \right\}$. 
\end{itemize}
In addition, we set $\Inf(m,d_A)=\underffp(m,d_A)$ for all $m\in\cX$, $A\subset \I$, $A\neq \emptyset$. Using Lemma \ref{lem:comput_env_inf}, we have that, for all $m,M\in\cX$, $a\in\cE$, $\Inf(m,a)\leq \FFppinf\left([m,M],a\right)$. Let $\tau^{(p)}_{\inf}$ be the time for $\Inf$ to hit $C_p$ (starting from $0$), then we have that $\tau^{(p)} \leq \tau^{(p)}_{\inf}$. Moreover, we can bound the mean of $\tau^{(p)}_{\inf}$ by:
$$ \Eb\left[\tau^{(p)}_{\inf}\right]\leq \textstyle{\frac{\Lambda}{\delta_p}} C_p.
$$
The arguments to show this are those of Lemma \ref{lem:tau} (\ref{subsubs:POS_slightly_general}) or Proposition \ref{prop:POS_zero} (Appendix \ref{sec:AppB}): we apply Theorem \ref{Foster} (\ref{subsubs:POS_slightly_general}) with different parameters, noting that $\Inf$ is also monotone. 
\endproof

\proof[Proof of Theorem \ref{th:algo_B}.]
Due to condition (iii), the projection of $H$ on $\I_0$ is a Markov chain. Condition (i) allows to apply Lemma \ref{lem:Foster_sup} to this Markov chain. 
Condition (ii) allows to apply Lemma \ref{lem:inf_capacity}, and the fact that $\I = \I_0 \cup \I_C$ concludes the proof.
\endproof



\section*{Acknowledgements}
This research is supported by the French National Research
Agency grant ANR-12-MONU-0019.



\bibliographystyle{apalike}




\bibliography{perfectATO} 

\end{document}